\documentclass[11pt,a4paper]{article}
\usepackage{mathptmx}
\usepackage[a4paper,margin=24mm]{geometry}
\usepackage{amsmath,amssymb,amsthm,mathtools,booktabs,tabularx,array,xcolor,graphicx}
\newcolumntype{Q}{>{\centering\arraybackslash}X}

\definecolor{primary}{RGB}{22, 66, 118}
\definecolor{accent}{RGB}{165, 55, 38}
\definecolor{softbg}{RGB}{244, 247, 251}

\usepackage[unicode,colorlinks=true,linkcolor=primary,citecolor=primary,urlcolor=primary]{hyperref}
\hypersetup{
  pdftitle={A Constructive Lower Bound for the Growth Constant of Ternary Cap Sets},
  pdfsubject={Combinatorics},
  pdfauthor={ZeHua Lu}
}

\newcommand{\F}{\mathbb{F}}

\newcommand{\Lang}{\mathcal{W}}
\DeclareMathOperator{\symdiff}{\triangle}

\theoremstyle{plain}
\newtheorem{theorem}{Theorem}[section]
\newtheorem{lemma}[theorem]{Lemma}
\newtheorem{proposition}[theorem]{Proposition}

\theoremstyle{definition}
\newtheorem{definition}[theorem]{Definition}
\newtheorem{example}[theorem]{Example}

\theoremstyle{remark}

\numberwithin{equation}{section}

\begin{document}

\begin{center}
  {\LARGE\bfseries\color{primary} A Constructive Lower Bound for the Growth Constant of Ternary Cap Sets\par}
  \vspace{0.55em}
  {\large Two-Color Recursion and Closed-Walk Synchronization\par}
  \vspace{0.65em}
  {\large ZeHua Lu\par}
  {\normalsize Independent Researcher\par}
  {\small\href{https://orcid.org/0009-0000-2307-5599}{ORCID 0009-0000-2307-5599}\par}
\end{center}
\vspace{0.6em}

\begin{abstract}
Let $r_3(d)$ denote the maximum size of a cap set in $\F_3^d$. We prove that the exponential growth constant $\Lambda = \lim_{d\to\infty} r_3(d)^{1/d}$ is strictly greater than $2.22920646813$. The main advance comes from combining two-color recursive synchronization with cross-layer joint filling of a seven-dimensional local alphabet that contains collinear triples (lines). By encoding marked positions as $14$-step closed walks on a two-state directed graph, any zero-sum triple with a fixed outermost color is forced to choose the identical closed walk—and hence the identical record—at every recursive layer, at a net cost of $1000$ dimensions per level. Because the open seven-dimensional coordinates are deferred to the outermost layer, all positions of the same type across all layers can be pooled per record and filled jointly via grouped binary-flag nodes. A $2\times 2$ transfer matrix governs the weighted record count, and spectral dominance over two tilted tail matrices converts this count into a capacity lower bound. Feeding each resulting bound back into the filling parameters iteratively amplifies the growth rate, while a companion-path score exclusion test further enlarges the admissible walk family. The final bound is certified by finite enumeration and exact integer comparisons.
\end{abstract}


\section{Introduction and Main Result}

In the finite vector space $\F_3^d$, three distinct points $x,y,z\in\F_3^d$ are collinear (i.e., form a three-term arithmetic progression) if and only if $x+y+z=0$. A subset $C\subseteq\F_3^d$ is called a \emph{cap set} if every triple $(x,y,z)\in C^3$ satisfying $x+y+z=0$ must have $x=y=z$. Let $r_3(d)$ denote the maximum cardinality of a cap set in $\F_3^d$. Because the Cartesian product of two cap sets is again a cap set, the sequence $\{r_3(d)\}_{d\ge 1}$ is supermultiplicative: $r_3(d+e)\ge r_3(d)r_3(e)$. By Fekete's superadditive lemma, the exponential growth constant (or capacity) of ternary cap sets exists and equals its supremum:
\[
  \Lambda = \lim_{d\to\infty} r_3(d)^{1/d} = \sup_{d\ge 1} r_3(d)^{1/d}.
\]

Edel's generalized product construction combines low-dimensional cap sets into higher dimensions \cite{edel}. Building on this framework, Tyrrell's method of extendable triples yielded a lower bound of $2.2180$ \cite{tyrrell}. Using an admissible set of $237\,984$ words found via program search, Romera-Paredes et al.\ obtained $2.2202$ \cite{funsearch}, and Zhai et al.\ pushed this to $2.2203$ with longer admissible words \cite{zhai}. Naslund subsequently announced a constructive lower bound of $2.22175$ \cite{naslund}. From the opposite direction, Ellenberg and Gijswijt used the polynomial method to prove $\Lambda \le \min_{t>0}(1+t+t^2)t^{-2/3} = 2.7551046\ldots$ \cite{eg}.

Hametner and Tyrrell recently proved that the per-dimension growth rate of any fixed finite cap set can be strictly improved in higher dimensions \cite{ht}. Obtaining a concrete, competitive exponential lower bound from such an improvement requires controlling both the coordinate cost of each lifting step and the number of points gained. In this paper, we establish the following bound.

\begin{theorem}\label{thm:main}
The growth constant of ternary cap sets satisfies
\[
  \Lambda > 2.229206468135675727981404.
\]
In particular, $r_3(d) > (2.22920646813)^d$ for all sufficiently large $d$.
\end{theorem}

\subsection*{Overview of the Construction and Sources of Gain}

Traditional product and admissible-set constructions require the local alphabet itself to be cap-free, which limits local encoding density. We bypass this bottleneck by coloring each coordinate position with a \emph{marked color} $T$ or an \emph{ordinary color} $O$ so that two copies of $T$ and one copy of $O$ at the same position can never sum to zero, while allowing the seven-dimensional local alphabets to contain collinear triples. The proof rests on two primary mechanisms, reinforced by recursive feedback and a companion-path refinement:

\begin{enumerate}
  \item \textbf{Closed-walk synchronization and single-slot recursion.} Section~\ref{sec:slots} formalizes \emph{two-color slots} with open coordinates, proving that zero-sum triples recover identical records under gluing and scalar deletion. Starting from a $3$-sunflower-free family $L$ of $352$ sets on $14$ points \cite{degkm,tableL}, Section~\ref{sec:walk} designates two bits as ports between consecutive rows. Any three non-identical compatible closed walks are locked into one of two mixed port states (Lemma~\ref{lem:ports}), both of which are eliminated by a conserved integer row score. Section~\ref{sec:rec} then embeds the previous-level construction into a single $7$-dimensional slot of the walk and deletes an output scalar from the first step, adding $1000$ net dimensions per layer. The resulting record counts are tracked by a $2\times 2$ transfer matrix $K(X,Y)$ (Theorem~\ref{thm:spectral}).

  \item \textbf{Cross-layer joint filling of line-containing alphabets.} Because open $7$-dimensional coordinates are deferred to the outermost layer, Section~\ref{sec:physical-filling} pools all open positions of the same type across all recursive layers for each synchronized record and fills them jointly using grouped binary-flag nodes (Lemma~\ref{lem:physical-node}) and two mixed grouping chains (Proposition~\ref{prop:filling}). Even when restricted to the base two-sided score intervals $\Lang_0$ and the trivial seed $z=2$, this cross-layer joint filling already lifts the capacity lower bound to $\Lambda > 2.228868987765\ldots$.

  \item \textbf{Recursive feedback and companion-path refinement.} Since the filling weights $(X,Y)$ depend monotonically on any certified lower bound $z < \Lambda$ (via the cap-set translate density $\beta=z/3$ and the zero-fiber rate $4+4z$), substituting each newly certified bound back into Proposition~\ref{prop:filling} and Theorem~\ref{thm:spectral} yields a second substantial jump, converging to $\Lambda > 2.229206468135\ldots$ over nine iterations (Appendix~\ref{sec:certificate}). In addition, tracking the achievable scores of two companion paths (Section~\ref{sec:walk}) enlarges $\Lang_0$ to $\Lang$ while preserving synchronization (Proposition~\ref{prop:sync}); this refinement provides a strict gain (already raising the $z=2$ first-round bound to $2.228868987790\ldots$) within the same finite integer certificate.
\end{enumerate}

\begin{figure}[!htbp]
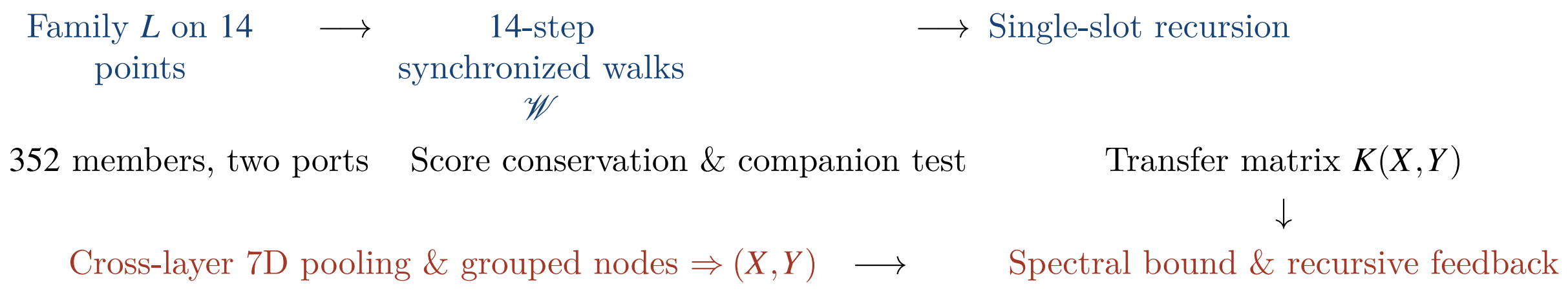

\centering
{\small
\renewcommand{\arraystretch}{1.25}
\begin{tabularx}{\linewidth}{@{}Q@{\enspace$\longrightarrow$\enspace}Q@{\enspace$\longrightarrow$\enspace}Q@{}}
\textcolor{primary}{Family $L$ on $14$ points} &
\textcolor{primary}{$14$-step synchronized walks $\Lang$} &
\textcolor{primary}{Single-slot recursion} \\
\multicolumn{1}{c}{$352$ members, two ports} &
\multicolumn{1}{c}{Score conservation \& companion test} &
\multicolumn{1}{c}{Transfer matrix $K(X,Y)$} \\
\multicolumn{2}{c}{} & \multicolumn{1}{c}{$\downarrow$} \\
\multicolumn{2}{c}{\textcolor{accent}{Cross-layer 7D pooling \& grouped nodes $\Rightarrow (X,Y)$}\quad$\longrightarrow$} &
\multicolumn{1}{c}{\textcolor{accent}{Spectral bound \& recursive feedback}}
\end{tabularx}}
\caption{Synchronized closed walks and single-slot recursion generate the record weight polynomial, while cross-layer seven-dimensional pooling supplies $(X,Y)$. Iteratively feeding the certified bound back into the filling weights drives the final bound.}
\label{fig:architecture}
\end{figure}

\section{Propagating Zero-Sum Conditions via Two-Color Slots}\label{sec:slots}

In this section we define two-color realizations with open coordinates and prove how zero-sum relations on physical vectors recover a common record under gluing and coordinate deletion.

Before stating the general definition, consider two low-dimensional prototypes. In $\F_3$, set $T=\{0\}$ and $O=\{-1,1\}$. Clearly, two copies of $T$ and one copy of $O$ cannot sum to zero, while three same-color scalars summing to zero must be equal. In $\F_3^3$, requiring exactly two positions to be $T$ and one to be $O$ yields the six points $\{\pm e_1,\pm e_2,\pm e_3\}$. Selecting three points with distinct color layouts produces two $T$'s and one $O$ in every coordinate; if two layouts match and the third differs, such a coordinate still exists. Thus any zero-sum triple must first agree in color layout, after which the one-dimensional property forces all scalar values to coincide. We now generalize this two-step recovery argument to realizations carrying unfilled coordinates.

\begin{definition}[Two-color slot]\label{def:two-color-slot}
Let $V$ be a fixed $\F_3$-vector space. A \emph{two-color slot} on $V$ consists of two classes of colored \emph{realizations}, labeled ordinary ($O$) and marked ($T$). Each realization is specified by three components: its color $c\in\{O,T\}$, a \emph{record} (listing the data already fixed at the slot, the types of open coordinates, and their ordering), and the values assigned to those \emph{open coordinates}; together these determine a physical vector in $V$. When combining multiple slots, distinct slots occupy pairwise disjoint coordinate blocks. The slot is required to satisfy two axioms:
\begin{enumerate}
  \item \textbf{Color exclusion:} If three realizations have physical vectors summing to zero, their colors cannot consist of exactly two $T$'s and one $O$.
  \item \textbf{Record and open-coordinate recovery:} If three realizations of the \emph{same color} have physical vectors summing to zero, they must share the same record, and at every corresponding open coordinate their three values sum to zero.
\end{enumerate}
\end{definition}

Note that triples in Definition~\ref{def:two-color-slot} may repeat realizations. This immediately yields disjointness and uniqueness.

\begin{lemma}[Color disjointness and record uniqueness]\label{lem:slot-unique}
In any two-color slot, the sets of physical vectors realized with color $O$ and color $T$ are disjoint. Moreover, within a fixed color class, each physical vector corresponds to at most one record and one assignment to its open coordinates.
\end{lemma}
\begin{proof}
If a physical vector $v\in V$ had both a $T$-realization and an $O$-realization, taking two copies of the $T$-realization and one of the $O$-realization would give $v+v+v=3v=0$, violating Axiom~1. Similarly, if two realizations $a,b$ of the same color have the same physical vector $v$, applying Axiom~2 to $(a,a,b)$ shows they share the same record, and at each open coordinate their values $(x,x,y)$ satisfy $x+x+y=0$, i.e., $y=x$. Hence $a=b$.
\end{proof}

We now show how to glue two-color slots along a synchronizing word family and delete an auxiliary scalar coordinate to form a higher-level two-color slot.

\begin{proposition}[Gluing and coordinate deletion]\label{prop:glue}
Let $\mathcal{W}$ be a family of words over $\{O,T\}$ on a fixed set of real positions, satisfying the property that any three words in $\mathcal{W}$ (with repetition allowed) having no real position with exactly two $T$'s and one $O$ must be identical. Place a two-color slot at each real position and take the direct sum over all real positions:
\begin{enumerate}
  \item \textbf{Global gluing synchronization:} If three glued vectors have physical sum zero, then they arise from the same word in $\mathcal{W}$, share identical subrecords at every real position, and have coordinatewise zero sum on all open coordinates.
  \item \textbf{Deleting a scalar coordinate:} Suppose the glued block contains a one-dimensional scalar coordinate whose zero or active status is determined by the record, and whenever it is active both signs $\pm 1$ can be realized while keeping the word and all other physical coordinates unchanged. Retain the realizations with scalar $0$ as the new ordinary color $O'$, retain the active realizations with scalar fixed to $+1$ as the new marked color $T'$, and delete this scalar coordinate. The resulting two classes of realizations in the remaining space form a two-color slot.
\end{enumerate}
\end{proposition}
\begin{proof}
For part~1, suppose three glued vectors sum to zero. Because distinct real positions occupy disjoint coordinate blocks, each real position sums to zero individually. Applying Axiom~1 at every real position rules out two $T$'s and one $O$ everywhere, so the hypothesis on $\mathcal{W}$ forces all three vectors to come from the same word. Thus all three realizations have the same color at every real position, and Axiom~2 recovers identical subrecords and coordinatewise zero sums on all open coordinates.

For part~2, consider three realizations whose projections after deleting the scalar coordinate sum to zero:
\begin{itemize}
  \item \emph{Axiom~1:} Suppose two have color $T'$ and one has color $O'$. Restore $+1$ and $-1$ to the deleted scalar of the first two realizations, respectively, and $0$ to the third. Because an active scalar can flip sign while preserving the word and all other physical coordinates, the first two lift to valid realizations in the original glued space. Since $(+1)+(-1)+0=0$, the lifted vectors sum to zero. By part~1, they recover the same word and subrecords, hence the same active flag on the deleted scalar—contradicting that two are active and the third is zero.
  \item \emph{Axiom~2:} Suppose all three have color $O'$ (respectively, $T'$). Restore $0$ (respectively, $+1$) to the deleted scalar of all three. Since $0+0+0=0$ and $1+1+1=0$ in $\F_3$, the lifted vectors sum to zero in the glued space. Part~1 again yields identical records and coordinatewise zero sums on all remaining open coordinates.
\end{itemize}
Hence the resulting structure satisfies both axioms of a two-color slot.
\end{proof}

In the seven-dimensional alphabets of Section~\ref{sec:physical-filling}, the position and value of each ordinary isolated point are written into the record, while the values of marked letters and active ordinary letters are left open for filling. Fixing the outermost color, we collect all open positions across all layers for each record in their original coordinate order, and choose a single cap set in the ambient space of those coordinates whose points lie in the Cartesian product of allowed values at each open position. The union of these filled point sets over all records remains a cap set: any zero-sum triple first recovers the same record layer by layer, so its three open-coordinate vectors sum to zero coordinatewise, and the cap-set property for that record forces the three vectors—and hence the three physical points—to be equal.

\section{Synchronizing 14-Step Closed Walks}\label{sec:walk}

This section constructs the synchronizing word family $\Lang$ required by Proposition~\ref{prop:glue}. Starting from a $3$-sunflower-free family on a $14$-element ground set, we use two coordinates as ports between consecutive rows. Closure forces any non-identical compatible triple into one of two port states, which are then eliminated by a conserved row score.

\subsection{Two Possible Port States}

Let $L$ be the family of $352$ sets on a $14$-element ground set constructed by Deuber, Erd\H{o}s, Gunderson, Kostochka, and Meyer \cite[Lemma~4.3]{degkm}; their two-port deletion gluing appears in \cite[Lemma~4.4]{degkm}. Following the public integer-mask transcription \cite{tableL}, we label the ground set $0,\ldots,13$ and identify each set with its binary characteristic row $r=(r_0,\ldots,r_{13})\in\{0,1\}^{14}$. For any three rows in $L$ that are not all identical, there exists at least one coordinate—called a \emph{witness coordinate}—containing exactly two $1$'s. Indeed, if the three rows are pairwise distinct, the $3$-sunflower-free property implies that the three pairwise intersections are not all equal, and any element in the difference of two intersections serves as a witness; if the three rows are $A,A,B$ with $A\ne B$, the antichain property ensures $A\setminus B \ne \varnothing$, any element of which is a witness.

Designating bits $0$ and $6$ as ports, define the \emph{source port} $i(r)$ and \emph{target port} $j(r)$ of row $r\in L$ by
\[
  i(r)=r_0, \qquad j(r)=1-r_6.
\]
We view each row $r\in L$ as a directed edge from $i(r)\in\{0,1\}$ to $j(r)\in\{0,1\}$. The remaining twelve bits $\{1,\ldots,5,7,\ldots,13\}$ are \emph{real slots}, where $1$ specifies $T$ and $0$ specifies $O$. A sequence of fourteen rows $r^{(0)},r^{(1)},\ldots,r^{(13)}\in L$ satisfying $j(r^{(t)})=i(r^{(t+1)})$ (indices modulo $14$) is called a \emph{closed walk} (or \emph{closed word}). Time $0$ is marked, so cyclic shifts generally produce distinct words. Three closed walks (or row words) are \emph{compatible} if no real slot contains exactly two $T$'s and one $O$.

\begin{lemma}[Port-state conservation]\label{lem:ports}
Suppose three $14$-step closed walks are compatible. If they are not all identical, then the number of $1$'s among the three source ports at time $t$ equals the same constant $s\in\{1,2\}$ for every $t\in\mathbb{Z}_{14}$.
\end{lemma}
\begin{proof}
Let $s_t\in\{0,1,2,3\}$ be the number of $1$'s among the three source ports at time $t$ (indices modulo $14$). Whenever the three rows at step $t$ are not all identical, compatibility rules out all twelve real slots as witness coordinates, so the witness must be either port bit $r_0$ or port bit $r_6$:
\begin{itemize}
  \item If $r_0$ is the witness, then $s_t=2$;
  \item If $r_6$ is the witness, then two rows have $r_6=1$; since $j(r)=1-r_6$ and $i(r^{(t+1)})=j(r^{(t)})$, this yields $s_{t+1}=1$.
\end{itemize}
Thus, whenever step $t$ has non-identical rows, either $s_t=2$ or $s_{t+1}=1$.

If $s_t=2$ at some time $t$, then the three target ports at time $t-1$ are mixed, so the rows at time $t-1$ are not all identical. Since $s_t=2\ne 1$, the implication above forces $s_{t-1}=2$. Propagating backward around the cycle gives $s_t=2$ for all $t$. Otherwise, $s_t\ne 2$ for all $t$. Picking any step $t$ where the three rows are not all identical forces $s_{t+1}=1$; hence the source ports at $t+1$ are mixed, so the rows at $t+1$ are also non-identical, forcing $s_{t+2}=1$. Propagating forward around the cycle gives $s_t=1$ for all $t$.
\end{proof}

\subsection{Score Intervals from the First Two Rows}

For each state $s\in\{1,2\}$, we enumerate all unordered one-step compatible triples of rows (allowing repetition) whose three source ports and three target ports each contain exactly $s$ ones. Both states yield $1350$ triples involving the exact same $330$ rows. The remaining $22$ rows are called \emph{blocking rows}; by Lemma~\ref{lem:ports}, any closed walk containing a blocking row cannot appear in a compatible triple of non-identical closed walks.

On the $330$ non-blocking rows, there exists an integer score function $e: L\setminus\{\text{blocking rows}\} \to \{0,1,2\}$ such that every compatible triple $(a,b,c)$ in state $s$ satisfies the conservation law
\begin{equation}\label{eq:step-conservation}
  e(a)+e(b)+e(c)=
  \begin{cases}
    4, & s=1,\\
    2, & s=2.
  \end{cases}
\end{equation}
The score assignment is supplied with the certificate. To enumerate compatible triples, one fixes the supports $A,B$ of two rows on the real slots; the third row's support must contain $A\cap B$ and avoid $A\symdiff B$, with remaining real slots arbitrary, and is then checked against ports and $L$.

For each non-blocking row, we record its source port, target port, number of $T$'s among the ten $7$-dimensional slots, number of active slots among the two scalar bits, active flag at bit $12$, and score $e$, yielding $18$ row types. Numbered lexicographically, pairs of consecutive initial rows $(r^{(0)},r^{(1)})$ form $162$ ordered \emph{head types}. Pairing consecutive compatible triples produces $468$ unordered compatible head-type triples in each state $s\in\{1,2\}$.

For each head type $h$, take the integer thresholds $\ell_h,u_h$ from the certificate. For a closed walk without blocking rows, let $G=\sum_{t=2}^{13} e(r^{(t)})$ be the score of its twelve-step suffix. Define the base family $\Lang_0$ to consist of all closed walks containing at least one blocking row, together with all blocking-free closed walks satisfying
\[
  \ell_h \le G \le u_h.
\]
The certificate thresholds satisfy
\begin{equation}\label{eq:head}
  u_h+u_j+u_k < 48 \quad (s=1), \qquad
  \ell_h+\ell_j+\ell_k > 24 \quad (s=2),
\end{equation}
for every compatible head-type triple $(h,j,k)$ in the respective state. Since three compatible blocking-free walks in state $s$ must have suffix scores summing to $12\times 4 = 48$ ($s=1$) or $12\times 2 = 24$ ($s=2$), inequalities \eqref{eq:head} and Lemma~\ref{lem:ports} imply that any three compatible walks in $\Lang_0$ are identical.

\subsection{Widening the Interval via Companion Paths}

For blocking-free walks outside $[\ell_h,u_h]$, we can admit a candidate walk $P$ by tracking the achievable suffix scores of the other two walks in any potential compatible triple—called its \emph{companion paths}.

Fix state $s\in\{1,2\}$. For a row $r$ in the suffix of $P$, define the $2\times 2$ set-valued matrix $B_s(r)$ whose entry $B_s(r)_{ab}$ ($a,b\in\{0,1\}$) is the set of scores $e(r')$ of all first-companion rows $r'$ with ports $a\to b$ that appear with $r$ in some allowed real one-step triple of state $s$. Compose two tables via set-sum matrix multiplication
\[
  (B*C)_{ab}=\bigcup_{v=0}^1\bigl(B_{av}+C_{vb}\bigr),
\]
where $U+V=\{u+v : u\in U,\, v\in V\}$. Multiplying the twelve suffix matrices gives
\[
  B_s(P)=B_s(r^{(2)})*B_s(r^{(3)})*\cdots*B_s(r^{(13)}),
\]
which contains the suffix score of every valid first companion suffix.

Let $P$ have head type $h$ and suffix score $G$, and set $D_1=48$, $D_2=24$. Enumerate all ordered compatible companion head-type pairs $(j,k)$ for $h$ in state $s$. The first companion suffix must return from its second row's target port $a$ to its first row's source port $b$. For each $w\in B_s(P)_{ab}\cap[0,24]$ with $D_s-G-w\in[0,24]$, the second companion suffix must score $D_s-G-w$. For both companion paths to obey their respective one-sided bounds from $\Lang_0$, we must simultaneously have
\begin{align}
  w \le u_j, \qquad D_1-G-w \le u_k &&(s=1),\label{eq:high}\\
  w \ge \ell_j, \qquad D_2-G-w \ge \ell_k &&(s=2).\label{eq:low}
\end{align}
If no compatible pair $(j,k)$ and entry $w\in B_s(P)_{ab}$ with $w, D_s-G-w\in[0,24]$ satisfy \eqref{eq:high} (for $s=1$) or \eqref{eq:low} (for $s=2$), we say $P$ \emph{passes the exclusion test for state $s$}. Because $B_s(P)$ allows the first companion to choose a local partner independently at each step, this relaxation safely covers all genuine compatible triples.

For each $s\in\{1,2\}$, choose a head-type subset $I_s$ such that every compatible head-type triple in state $s$ contains at most one element of $I_s$ counted with multiplicity. We enlarge $\Lang_0$ to $\Lang$ by adding all blocking-free walks $P$ satisfying either
\[
\begin{array}{ccl}
  G > u_h & \text{and } h\in I_1, & P\text{ passes the state-$1$ exclusion test},\\[2pt]
  G < \ell_h & \text{and } h\in I_2, & P\text{ passes the state-$2$ exclusion test}.
\end{array}
\]
The sets $I_1,I_2$ are chosen independently: each head type is weighted by the sum of main-matrix weights of its newly admitted walks across all four parent-child color combinations; we form the conflict graph whose edges connect head types co-occurring in a compatible triple, remove self-conflicting head types, and greedily select vertices while deleting their neighbors. The certificate checks the multiplicity condition against all $468$ triples in each state after every selection.

\begin{proposition}[Synchronization of $\Lang$]\label{prop:sync}
Any three compatible closed words in $\Lang$ are identical.
\end{proposition}
\begin{proof}
Suppose for contradiction that $\Lang$ contains three compatible walks that are not all identical. By Lemma~\ref{lem:ports}, they lie in state $s=1$ or $s=2$, and hence none contains a blocking row.

If $s=1$, their suffix scores sum to $D_1=48$. By \eqref{eq:head}, at least one walk $P$ (of head type $h$) has $G > u_h$. By construction of $\Lang$, we have $h\in I_1$ and $P$ passes the state-$1$ exclusion test. By the multiplicity condition on $I_1$, the other two walks have head types $j,k\notin I_1$, so neither can exceed its upper bound (even if admitted from below via $I_2$). Their actual suffix scores therefore give a genuine solution to \eqref{eq:high}, contradicting that $P$ passed the state-$1$ exclusion test.

If $s=2$, the suffix sum $D_2=24$ and \eqref{eq:head} analogously force one walk to satisfy $G < \ell_h$ with $h\in I_2$, while the remaining two obey their lower bounds and yield a genuine solution to \eqref{eq:low}, again a contradiction.
\end{proof}

\noindent\begin{minipage}{\linewidth}\begin{example}
Writing $14$-bit masks in decimal, consider the closed walk
\[
  (12302,6769,5734^{(7)},10652,2463^{(3)},5731),
\]
where parenthesized superscripts denote consecutive repetitions. This word has head type $45$ and suffix score $G=16$, exceeding $u_{45}=15$. Fixing $P$, the relaxed first-companion score table $B_1(P)$ and the scores required by \eqref{eq:high} across all compatible ordered companion head types are:
\[
\begin{array}{ccl}
 \text{Suffix ports} & B_1(P)_{ab} & \text{Scores satisfying both upper bounds}\\
 \hline
 0\to0 & \{10,12,\ldots,24\} & \{15\}\\
 0\to1 & \{9,11,\ldots,23\}  & \{16\}\\
 1\to0 & \{9,11,\ldots,23\}  & \{16\}\\
 1\to1 & \{8,10,\ldots,22\}  & \{17\}
\end{array}
\]
In each row, the last two columns have opposite parity and hence empty intersection, so $P$ passes the state-$1$ exclusion test. In both Round~$0$ and Round~$9$ of Appendix~\ref{sec:certificate}, the independent set $I_1$ chosen by the main matrix weights contains head type $45$, so $P\in\Lang\setminus\Lang_0$ in both rounds. As a cross-check, tracking two genuinely closed companion paths simultaneously along the raw rows yields $64$ head-type and score combinations, none of which satisfies both upper bounds.
\end{example}\end{minipage}\par

\section{Single-Slot Recursion and Growth Rate}\label{sec:rec}

We now use $\Lang$ to recursively construct two-color slots, count records via a $2\times 2$ transfer matrix, and deduce a capacity lower bound from spectral radii.

\subsection{Recursive Architecture and Dimension Accounting}

Each row of $\Lang$ has twelve real slots: we assign bits $12$ and $13$ to one-dimensional scalars (ordinary or marked) and each of the remaining ten real slots to a seven-dimensional alphabet, occupying
\[
  10\times 7 + 2\times 1 = 72 \text{ dimensions per row}.
\]
Section~\ref{sec:physical-filling} constructs two seven-dimensional alphabets in $\F_3^7$ such that:
\begin{itemize}
  \item At an \emph{ordinary} $7$-dimensional slot, $64$ isolated points are recorded directly, while the active ordinary letters are left open with formal weight $Y$, giving total formal weight $64+Y$.
  \item At a \emph{marked} $7$-dimensional slot, there are no isolated points, and the open filling is tracked by formal weight $X$.
\end{itemize}

The recursion starts from a single point in $\F_3^0=\{0\}$ colored $O$, with the initial $T$ class empty. At depth $n$, along a $14$-step closed walk (steps $t\in\{0,\ldots,13\}$), we perform two fixed operations:
\begin{enumerate}
  \item \textbf{Embedding the previous level:} Step~$7$, bit~$5$ (a $7$-dimensional slot) is reserved to house the entire two-color construction from depth $n-1$; its color ($T$ or $O$) selects the output color of the embedded level.
  \item \textbf{Deleting the output scalar:} Step~$0$, bit~$12$ serves as the output scalar. If the original row has $T$ at this bit, the scalar is $0$ and the new output color is $O$; if the original row has $O$, the scalar takes $\pm 1$, of which we retain $+1$ and assign the new output color $T$. This coordinate is then deleted as in Proposition~\ref{prop:glue}.
\end{enumerate}
By Propositions~\ref{prop:glue} and~\ref{prop:sync}, every depth yields a valid two-color slot. Fourteen rows total $14\times 72 = 1008$ coordinates; subtracting the $7$-dimensional embedded slot and the $1$-dimensional deleted scalar gives a net increase of
\[
  14\cdot 72 - 7 - 1 = 1000
\]
dimensions per level. All records of depth $n$ are thus realized in $\F_3^{1000n}$.

\subsection{Transfer Matrix and Spectral Bound}

Index rows and columns by $(O,T)$ and define the $2\times 2$ matrix $K(X,Y)$ whose $(a,b)$-entry is the sum of weights of all allowed closed walks with output color $a$ and embedded-slot color $b$. Every real slot other than the embedded slot and the deleted scalar contributes its factor to the product: $X$ for a marked $7$-dimensional slot, $64+Y$ for an ordinary $7$-dimensional slot, $2$ for an ordinary scalar, and $1$ for a marked scalar. Thus $(K(X,Y)^n)_{00}$ counts all $n$-level records whose innermost and outermost colors are both $O$, with the monomial $X^A Y^B$ recording the total numbers $A$ of marked and $B$ of active ordinary $7$-dimensional slots across the entire record.

All four entries of $K(X,Y)$ are strictly positive: taking mask $5734$ or $2302$ at step~$0$ selects output color $O$ or $T$, taking mask $2010$ or $2302$ at step~$7$ selects embedded color $O$ or $T$, and taking the blocking row $2302$ at all other steps forms valid $0\to 0$ self-loops in all four cases.

Proposition~\ref{prop:filling} applies to records whose counts $(A,B)$ lie in the closed convex cone \eqref{eq:cone} spanned by $(137,6)$ and $(9,718)$. We control records outside this cone via two tilted matrices.

\begin{theorem}[Spectral bound for single-slot recursion]\label{thm:spectral}
Let $X,Y>0$. Suppose there exists $c>0$ such that every record satisfying
\begin{equation}\label{eq:cone}
  9B \le 718A, \qquad 6A \le 137B
\end{equation}
admits a cap set of size at least $cX^A Y^B$ on its open positions. Fix $t>1$ and set
\[
  K_0 = K(X,Y), \qquad K_- = K(Xt^{-718}, Yt^9), \qquad K_+ = K(Xt^6, Yt^{-137}).
\]
If $\rho(K_0) > \max\bigl(\rho(K_-), \rho(K_+)\bigr)$, then
\[
  \Lambda \ge \rho(K_0)^{1/1000}.
\]
\end{theorem}
\begin{proof}
On any $n$-level record with $O$ at both the innermost and outermost levels and counts $(A,B)$, its weight in $(K_0^n)_{00}$ is $X^A Y^B$. If $9B > 718A$, the substitution in $K_-$ multiplies this weight by $t^{9B-718A} > 1$; if $6A > 137B$, the substitution in $K_+$ multiplies it by $t^{6A-137B} > 1$. Since $K_-$ and $K_+$ have nonnegative entries, the total weight of records satisfying the closed cone condition \eqref{eq:cone} is at least
\[
  (K_0^n)_{00} - (K_-^n)_{00} - (K_+^n)_{00}.
\]
Choosing the hypothesized cap set of size at least $cX^A Y^B$ for each such record, Propositions~\ref{prop:glue} and~\ref{prop:sync} ensure that these cap sets are pairwise disjoint and their union is a cap set in $\F_3^{1000n}$. Hence
\[
  r_3(1000n) \ge c\Bigl((K_0^n)_{00} - (K_-^n)_{00} - (K_+^n)_{00}\Bigr).
\]
By the Perron--Frobenius theorem applied to the strictly positive matrix $K_0$, $(K_0^n)_{00}/\rho(K_0)^n$ converges to a positive constant as $n\to\infty$, whereas both tail terms vanish when divided by $\rho(K_0)^n$. Taking $1000n$-th roots and letting $n\to\infty$ yields $\Lambda \ge \rho(K_0)^{1/1000}$.
\end{proof}

\section{Seven-Dimensional Alphabets and Path Filling}\label{sec:physical-filling}

This section constructs the ordinary alphabet $U$ and marked alphabet $M$ in $\F_3^7$ and fills all open seven-dimensional positions of a record simultaneously via grouped flag nodes.

First, for any $1 < z < \Lambda$, taking products of a finite cap set whose per-dimension root exceeds $z$ and padding remaining coordinates with $\{0,1\}^r$ gives a constant $L_z \ge 1$ depending only on $z$ such that for all $d\ge 0$,
\begin{equation}\label{eq:Lz-bound}
  r_3(d) \ge \frac{z^d}{L_z}.
\end{equation}

\subsection{Geometry of the Two Seven-Dimensional Alphabets}

Project the $8$-dimensional $512$-point cap set of Romera-Paredes et al.\ \cite[Figure~4 and Appendix~E.2]{funsearch} by deleting its zeroth coordinate. Exactly $160$ image points have two preimages; let $T_0 \subset \F_3^7$ be this set of $160$ points. Define
\[
  U = \F_3^7 \setminus (-T_0 - T_0), \qquad M = T_0 \cup E \cup C \cup D_+ \cup D_-,
\]
where $|E|=4$, $|C|=12$, and $|D_+|=|D_-|=3$ (explicit coordinates and coloring checks appear in Appendix~\ref{app:seven-alphabet}). Pairwise verification yields
\[
  |U|=200, \qquad |M|=182, \qquad M\cap U=\varnothing, \qquad (-M-M)\cap U=\varnothing.
\]
The identity $(-M-M)\cap U=\varnothing$ rules out zero-sum triples with two marked letters and one ordinary letter, establishing Axiom~1 of Definition~\ref{def:two-color-slot}.

Both $U$ and $M$ contain collinear triples (called \emph{proper lines}), structured into clean fiber hierarchies summarized in Table~\ref{tab:alphabet-structure}.

\begin{table}[!htbp]
\centering
\caption{Fiber decomposition of the ordinary alphabet $U$ and marked alphabet $M$ in $\F_3^7$}
\label{tab:alphabet-structure}
\begin{tabularx}{\linewidth}{@{}lllX@{}}
\toprule
Alphabet & Size & First-level partition & Second-level \& internal geometry \\
\midrule
Ordinary $U$ & $200$ & $64$ isolated $+$ $136$ active & Active points partition into $(64, 8, 64)$, each color a cap set \\
Marked $M$ & $182$ & Neg.\ $83$, Zero $16$, Pos.\ $83$ &
  Each nonzero fiber has $24$ isolated points and $59$ points partitioned into $(29,1,29)$.
  Each $29$-point fiber has $6$ isolated points and an $11$-spoke star.
  Zero fiber has $4$ isolated points ($E$) and $12$ points ($C$, $8$ lines). \\
\bottomrule
\end{tabularx}
\end{table}

\begin{itemize}
  \item \textbf{Structure of $U$:} $U$ has $64$ isolated points lying on no proper line; the remaining $136$ active points are colored $-1,0,+1$ with sizes $64,8,64$ so that every proper line has one point of each color, making each color class a cap set. Any zero-sum triple using an isolated point must repeat that point, so isolated points are recorded directly, leaving $128$ nonzero and $8$ zero choices at active ordinary positions.
  \item \textbf{Structure of $M$:} $M$ partitions into negative, zero, and positive fibers of sizes $83,16,83$ so that the three colors on every proper line sum to zero. Each nonzero fiber has $24$ isolated points, and the remaining $59$ points further partition into $(29,1,29)$ preserving zero-sum on lines. Each $29$-point second-level fiber consists of $6$ isolated points and an \emph{eleven-spoke star} ($11$ proper lines sharing a single center and otherwise pairwise disjoint, totaling $23$ points).
  \item \textbf{Capacity of the zero fiber $E\cup C$:} The $4$ points of $E$ are isolated, while the $8$ proper lines in $C$ each map onto $\F_3$ under $f(v)=v_1+v_6$, with each value in $\F_3$ having $4$ preimages in $C$. Pulling back a maximum $k$-dimensional cap set via $f$ yields a cap set of size $4^k r_3(k)$ on $k$ positions of $C$. Summing over all $\binom{d}{k}$ placements of $k$ positions in $C$ and $4^{d-k}$ choices on the remaining $d-k$ positions in $E$, and applying \eqref{eq:Lz-bound}, $d$ pooled zero-fiber positions accommodate at least
  \begin{equation}\label{eq:zero-fiber-cap}
    \sum_{k=0}^d \binom{d}{k} 4^{d-k} \cdot 4^k r_3(k)
    \ge \frac{1}{L_z}\sum_{k=0}^d \binom{d}{k} 4^{d-k} (4z)^k
    = \frac{(4+4z)^d}{L_z}
  \end{equation}
  points.
\end{itemize}

\subsection{Letter Weights and Fixed-Sum Grouped Nodes}

For the scalar alphabet $\F_3$ with parameter $x\ge 0$, assign weights $1/2, 1/2, x$ to $+1, -1, 0$. The weight of a word $w\in\F_3^d$ is the product of its coordinate weights, and the weight of a cap set is the sum of weights of its words. Define
\[
  F(x) = \sup_{d\ge 1} \Bigl(\max_{C\subseteq\F_3^d\text{ cap set}} \sum_{w\in C} \prod_{i=1}^d \mathrm{wt}_x(w_i)\Bigr)^{1/d}.
\]
By supermultiplicativity, for every $0 < s < F(x)$ there exists $c_s > 0$ such that every length $d$ admits a cap set of weight at least $c_s s^d$.

\begin{itemize}
  \item \textbf{Capacity of an $11$-spoke star:} Label the two endpoints of each spoke $+1$ and $-1$ and the center $0$. Pulling back a scalar cap set gives $11$ choices per nonzero symbol and $1$ choice for $0$, yielding a growth lower bound of $22F(1/22)$.
  \item \textbf{Merging general fibers:} Let $A$ be a finite alphabet with coloring $\phi:A\to\F_3$ satisfying $\phi(a)+\phi(b)+\phi(c)=0$ whenever $a+b+c=0$. Suppose the two nonzero fibers have positive capacity lower bounds $v_+, v_-$ and the zero fiber has positive capacity lower bound $Z$; set $V=v_+ + v_-$. Independently swap the $+1$ and $-1$ labels at each position with equal probability. Fix a scalar cap set and sum only over its symbol words; for each word, pool positions of the same physical fiber and take a cap set in each product space. Any physical zero-sum triple has zero-sum symbol words, forcing the symbol words to agree, after which the fiber cap sets force the physical words to agree. Taking expectations over sign flips gives average weight $V/2$ per nonzero symbol and $Z$ per zero symbol. Using strict lower bounds incurs one length-independent constant per fiber, yielding capacity lower bound $VF(Z/V)$; if an endpoint is attained, multiply the three bounds by $1-\varepsilon$ and let $\varepsilon\downarrow 0$.
\end{itemize}

\begin{lemma}[Fixed-sum grouped nodes]\label{lem:physical-node}
For each $b\in\{144,182,183,252,641\}$, there exists a randomized partial coloring $\chi:\F_3^b\dashrightarrow\F_3$ with acceptance probability lower bounds listed in Table~\ref{tab:node-probs}. It accepts only one-zero words (assigned nonzero colors) and all-nonzero words (assigned color $0$), and satisfies:
\begin{enumerate}
  \item \textbf{Equal-sum zero-sum recovery:} If three accepted input words sum to zero coordinatewise and share the same coordinate sum $\pi(w)=\sum_{i=1}^b w_i \in \F_3$, then their three output colors sum to zero; moreover, if the three output colors are equal, the three input words are identical.
  \item \textbf{Acceptance lower bounds:} When the remaining nonzero coordinates are chosen independently and uniformly from $\{\pm 1\}$, any one-zero word with a specified zero coordinate is accepted with probability at least $\eta_b$, and an all-nonzero word is accepted with probability at least $B_b$.
\end{enumerate}
\end{lemma}

\begin{table}[!htbp]
\centering
\caption{Acceptance probability lower bounds for grouped nodes (exact rationals scaled by $10^{12}$ and floored)}
\label{tab:node-probs}
\begin{tabular}{@{}rrr@{}}
\toprule
Input length $b$ & $\lfloor 10^{12}\eta_b\rfloor$ & $\lfloor 10^{12}B_b\rfloor$ \\
\midrule
$144$ & $313\,436\,531\,305$ & $328\,092\,514\,121$ \\
$182$ & $313\,436\,531\,305$ & $328\,035\,264\,355$ \\
$183$ & $313\,436\,531\,305$ & $328\,032\,524\,173$ \\
$252$ & $313\,436\,531\,305$ & $327\,984\,670\,844$ \\
$641$ & $313\,436\,531\,305$ & $327\,904\,358\,794$ \\
\bottomrule
\end{tabular}
\end{table}

The randomized binary-flag construction and probability derivations for Lemma~\ref{lem:physical-node} are given in Appendix~\ref{app:flag-node} (both bounds come from a single orientation choice, and computations use exact rationals).

When $g$ nodes are used in parallel at a grouping layer, we collect their input coordinate sums into a vector in $\F_3^g$ and require it to lie in a fixed random translate of a maximum cap set in $\F_3^g$. For any fixed record, uniform translation retains a fraction $r_3(g)/3^g \ge \beta^g/L_z$, where $\beta=z/3$. For three coordinatewise zero-sum records, the sum vectors are also zero-sum and hence equal by the cap-set property, so Lemma~\ref{lem:physical-node} propagates zero-sum outputs upward layer by layer. Restricting the top layer to nonzero outputs forces the three top-layer words to agree, and Lemma~\ref{lem:physical-node} then recovers identical input symbol words downward layer by layer. Translates across layers are drawn independently; for a fixed number of layers, each layer incurs at most one constant factor $L_z$, independent of $g$.

\subsection{Two Mixed Grouping Chains and the Filling Proposition}

Consider a $b$-input grouped node whose two nonzero symbols per input have total normalized weight $1$ and whose input zero weights sum to $t$. Because the node contributes one coordinate sum, it incurs one factor $\beta$ in the parallel translate estimate, giving output weight lower bounds $\beta\eta_b t$ (nonzero) and $\beta B_b$ (zero), with signs flipped independently before the next layer.

Chaining such nodes by feeding each node's output into the $b$-th input of the next node—with the remaining $b-1$ fresh inputs having zero-weight sum $t$—yields the transfer matrix and spectral radius
\[
  K_b(t) = \beta \begin{pmatrix} \eta_b t & \eta_b \\ B_b & 0 \end{pmatrix},
  \qquad H_b(t) = \rho\bigl(K_b(t)\bigr).
\]
Each chain is initialized with a uniform nonzero starting symbol fixed across all inputs of that chain, and only nonzero outputs are kept at the top. For $t>0$, $K_b(t)^2$ is strictly positive, so for any $0 < h < H_b(t)$ and sufficiently large chain length $d$, we have $(K_b(t)^d)_{11} > h^d$. Parallel chains of the same length share a single cap-set translate per height.

Choosing $s,v > 0$ such that
\[
  s^{251} < H_{252}(251/22), \qquad D = 12 + 44s, \qquad v^{640} < H_{641}(640/D),
\]
the $252$-input chain gives a global lower bound $F(1/22)\ge s$ for the $11$-spoke star, so the two second-level nonzero fibers ($6$ isolated points plus a star each) have combined growth at least $D=12+44s$; the $641$-input chain then gives a global lower bound $F(1/D)\ge v$. Including the $24$ isolated points in each first-level nonzero fiber of $M$, the combined capacity lower bound of the two first-level nonzero fibers is
\[
  N = 48 + 2Dv, \qquad x = \frac{4+4z}{N}, \qquad y = \frac{1}{16}.
\]
Factoring out $N^A$ and $128^B$ from $A$ marked and $B$ active ordinary $7$-dimensional positions leaves normalized scalar weights with zero weights $x$ and $y$, respectively.

\begin{proposition}[Seven-dimensional filling on a full record]\label{prop:filling}
Set
\[
\begin{aligned}
  H_1 &= H_{144}(137x + 6y),\\
  U_0 &= \beta\eta_{183}(3x + 180y), \qquad
  R_0 = \frac{B_{183}}{\eta_{183}(3x + 180y)},\\
  H_2 &= U_0^3 H_{182}(3R_0 + 178y).
\end{aligned}
\]
Choose $p,q>0$ satisfying
\[
  p^{137}q^6 < H_1, \qquad p^9 q^{718} < H_2,
\]
and set $X=Np$, $Y=128q$. For any $(A,B) \in \mathbb{R}_{\ge 0}(137,6) + \mathbb{R}_{\ge 0}(9,718)$, given the $A$ marked positions, $B$ active ordinary positions, and the ordinary isolated-point record, there exists a cap set of size at least $cX^A Y^B$ on the actual seven-dimensional alphabet, where $c>0$ is independent of $(A,B)$ and position ordering. Under outer closed-walk synchronization, each record may be filled independently, giving ordinary-position formal weight $64+Y = 64+128q$.
\end{proposition}
\begin{proof}
We build two mixed chains corresponding to the two extreme rays $(137,6)$ and $(9,718)$:
\begin{itemize}
  \item \textbf{Type-1 chain for $(137,6)$:} Each step feeds $137$ marked inputs, $6$ active ordinary inputs, and the previous node's output into a $b=144$ node, giving per-step growth at least $H_1$.
  \item \textbf{Type-2 chain for $(9,718)$:} At the bottom, three independent groups of $3$ marked and $180$ active ordinary inputs each feed into a $b=183$ node, with output weight lower bounds $(U_0, \beta B_{183}) = U_0(1, R_0)$. Their three outputs, $178$ fresh active ordinary inputs, and the previous main-chain output feed into a $b=182$ node, consuming $(9,718)$ seven-dimensional positions per step with growth at least $H_2$.
\end{itemize}
By the strict power inequalities, we choose a single finite chain length $d$ so that the top nonzero weights of both full chains exceed $(p^{137}q^6)^d$ and $(p^9 q^{718})^d$, respectively. Writing $(A,B) = u(137,6) + w(9,718)$ ($u,w\ge 0$) and rounding $u,w$ down to integer multiples of $d$ leaves fewer than $146d$ marked and $724d$ active ordinary leftover positions, which are fixed to nonzero symbols at the cost of a constant depending only on $d$.

Regardless of outer nesting depth, all seven-dimensional positions of the record are pooled once and filled together; a fixed number $m$ of grouping layers costs at most $L_z^m$. By linearity of expectation, some choice of flags, sign flips, and translates achieves the expected bound. Multiplying back $N^A 128^B$ and pooling identical physical fibers across all positions (incurring one length-independent constant per fiber) yields at least $cX^A Y^B$ points.

For the cap-set property, three zero-sum filled points first recover the same outer record and identical ordinary isolated points. On the remaining positions, the physical letter colorings sum to zero; the shared translates and Lemma~\ref{lem:physical-node} propagate zero sums upward to identical top-layer outputs and then recover identical input symbol words downward, placing all three points in the same physical fiber at every position. The fiber cap sets (using \eqref{eq:zero-fiber-cap} for the zero fiber) then force the three points to coincide.
\end{proof}

The function \texttt{physical\_prices} in \path{outer/verify_balanced_loop_quotient.py} computes exact rationals $N,p,q$ and verifies both strict inequalities.

\begin{proof}[Proof of Theorem~\ref{thm:main}]
The $3$-dimensional $9$-point cap set gives $9^{1/3}\le \Lambda$, so $z=2$ is a valid strict seed. Iteratively applying Proposition~\ref{prop:filling} and Theorem~\ref{thm:spectral}—verified by the integer certificate in Appendix~\ref{sec:certificate}—yields in Round~$9$ with $f=2.2292$ a normalized main matrix $K_0/f^{1000}$ satisfying
\[
  \rho(K_0/f^{1000}) \ge R = \frac{10029057595358243801843755621}{10^{28}},
\]
while both tilted tail matrices have spectral radius strictly below $R$, and
\[
  \Bigl(\frac{2.229206468135675727981404}{f}\Bigr)^{1000} < R.
\]
By Theorem~\ref{thm:spectral}, $\Lambda \ge f R^{1/1000} > 2.229206468135675727981404$, and the asymptotic statement for $r_3(d)$ follows from $\Lambda = \lim_{d\to\infty} r_3(d)^{1/d}$.
\end{proof}

\clearpage
\appendix

\section{Integer Certificate from the 9-Point Cap Set to Theorem~\ref{thm:main}}\label{sec:certificate}

This appendix details the seed construction, exact rational parameters, and pure integer comparisons used to certify Theorem~\ref{thm:main}. Throughout, we fix decimal grid precision $G=10^{24}$, matrix scale $S=10^{28}$, normalization factor $f=2.2292=5573/2500$, and tilt parameter $t=1000001/1000000$, using the same two-row score threshold table in every round.

\subsection{The 9-Point Seed and All-Dimension Lower Bounds}

In $\F_3^3$, the paraboloid $P=\{(a,b,a^2+b^2):a,b\in\F_3\}$ is a $9$-point cap set: any three collinear points have first two coordinates $u-h,u,u+h$, so their last coordinates sum to $2(h_1^2+h_2^2)$, which vanishes in $\F_3$ only when $h=0$. Products give $\Lambda \ge 9^{1/3} > 2$; since $\{0,1\}^j$ is a cap set for all $j\ge 0$, the seed $z=2$ satisfies $r_3(j)\ge 2^j$ with $L_2=1$.

In general, for any certified $1<z<\Lambda$, pick a cap set $D\subset\F_3^m$ with $|D|>z^m$. Writing $j=km+r$ ($0\le r<m$), the product $D^k\times\{0,1\}^r$ gives
\[
  r_3(j)\ge C_z z^j, \qquad C_z=\min_{0\le r<m}(2/z)^r > 0,
\]
ensuring \eqref{eq:Lz-bound} holds in every subsequent round.

\subsection{Rational input for one round}\label{appcert:prices}

A round takes a lower bound \(z<\Lambda\) and returns a new lower bound on
the grid \(G^{-1}\mathbb Z\). Its finite input is given in the electronic
supplement \path{appendix_a_data.json}. The array \texttt{rows} lists the
352 pairs (row mask, score), using \texttt{null} for a blocking row.
The arrays \texttt{lower} and \texttt{upper} give the 162 head thresholds
in the order specified below. The five pairs \((\eta_b,B_b)\) in
\texttt{node\_rates} are the exact rationals of Appendix~\ref{app:flag-node}.
For each seed, \texttt{rounds} gives the four positive rational witnesses
\(s,v,p,q\) used in Proposition~\ref{prop:filling}.

These witnesses can be checked without approximating a spectral radius.
For \(h>0\) the characteristic polynomial of \(K_b(u)\) has one positive
root and one negative root, so
\begin{equation}\label{appcert:price-test}
 h<H_b(u) \quad\Longleftrightarrow\quad
 h(h-\beta\eta_b u)<\beta^2\eta_b B_b,\qquad \beta=z/3 .
\end{equation}
Set \(D=12+44s,\ N=48+2Dv,\ x=(4+4z)/N,\ y=1/16\) and
\(U_0=\beta\eta_{183}(3x+180y),\ R_0=\beta B_{183}/U_0\).
The four applications of \eqref{appcert:price-test}, in order, are
\[
\begin{array}{c@{\qquad}l@{\qquad}l}
 b & u & h\\ \hline
 252 & 251/22 & s^{251}\\
 641 & 640/D & v^{640}\\
 144 & 137x+6y & p^{137}q^6\\
 182 & 3R_0+178y & p^9q^{718}/U_0^3 .
\end{array}
\]
They certify the filling prices \(X=Np\) and \(Y=128q\). For example,
in Round~9 the following integers give all four witnesses exactly.
\[
\begin{aligned}
 10^{32}s&=100393292983991794087700885356835, &
 10^{32}v&=100153847849986176176512064954925,\\
 10^{32}p&=100697311688384909000046865394478, &
 10^{32}q&=100541177178899363341387036079334 .
\end{aligned}
\]

\subsection{From one row to a closed-walk matrix}\label{appcert:walk}

\paragraph{Row types and compatible heads.}
Read mask \(r\) from its least significant bit, so that
\(r_d=\lfloor r/2^d\rfloor\bmod 2\). Write \(\bot\) for the blocking
symbol. With \(\mathcal D=\{1,2,3,4,5,7,8,9,10,11\}\), put
\[
 i_r=r_0,\quad j_r=1-r_6,\quad k_r=\sum_{d\in\mathcal D}r_d,\quad
 \nu_r=2-r_{12}-r_{13},\quad a_r=1-r_{12}.
\]
The type \(x(r)\) of a non-blocking row is the zero-based index of its
tuple \((i_r,j_r,k_r,\nu_r,a_r,e(r))\) among the 18 distinct tuples in
lexicographic order. Write \(i_x,j_x,e_x\) for the corresponding
components. A blocking row has type and score \(\bot\).
Index the pairs \((x,y)\) with \(j_x=i_y\) in lexicographic order,
writing \(h(x,y)\in\{0,\ldots,161\}\) for the head index.

To reconstruct the ordered row relation in state \(\sigma\in\{1,2\}\),
enumerate unordered triples of rows with repetition. Retain a triple
when its source bits sum to \(\sigma\), its target bits sum to \(\sigma\),
and each real slot has \(0,1\), or \(3\) ones. Insert all distinct
permutations. Before permutation there are 1350 triples in each state,
using precisely the 330 scored rows. Check the score sum in
\eqref{eq:step-conservation} on each triple. For each row \(r\) collect
the ordered pairs \((x(r'),x(r''))\) that occur with \(r\) first.
The resulting sets agree for all rows of the same type \(x\) in each
state. This finite equality justifies using their common set \(C_\sigma(x)\)
at both steps of the head.

For \(h=h(x,y)\), take \((c,d)\in C_\sigma(x)\) and
\((v,w)\in C_\sigma(y)\) with \(j_c=i_v,\ j_d=i_w\). Record
\[
 (h',h'',\alpha,\gamma)=\bigl(h(c,v),h(d,w),j_v,i_c\bigr)
 \quad\text{in }\mathcal H_\sigma(h).
\]
The two port equalities are precisely the joins needed for the companion
rows, so this enumeration covers every compatible two-step head. The
first companion must return from port \(\alpha\) to port \(\gamma\)
during the suffix. Sorting each triple \((h,h',h'')\) and
removing duplicates gives 468 triples per state. On these triples check
\eqref{eq:head}, counting a repeated head at each occurrence. Also check
\(0\le\ell_h\le u_h\le24\) for each \(h\).

\paragraph{Row weights and initial rounding.}
For a single round, run the recursion three times with the following
marked and ordinary prices and integer rounding maps.
\[
\begin{array}{c@{\qquad}c@{\qquad}c@{\qquad}c}
 & m & o & \mathcal R(\xi)\\ \hline
 K_0 & X & 64+Y & \lfloor \xi\rfloor\\
 K_- & Xt^{-718} & 64+Yt^9 & \lceil \xi\rceil\\
 K_+ & Xt^6 & 64+Yt^{-137} & \lceil \xi\rceil
\end{array}
\]
The 64 isolated ordinary choices carry no tilt. Fix an output color
\(a\) and an embedded color \(b\), in the order \(0=O,\ 1=T\).
At time \(d=0\) require \(a_r=a\), and at \(d=7\) require \(r_5=b\).
Let \(\delta=\mathbf 1_{\{d=7\}}\) and \(k=k_r-\delta b\).
The weight of an eligible row is
\begin{equation}\label{appcert:row-weight}
 W_d(r)=
 \frac{2^{\,\nu_r-a\mathbf 1_{\{d=0\}}}\,m^{k}o^{\,10-\delta-k}}
      {f^{\,7(10-\delta)+2}} .
\end{equation}
At time zero only the sign \(+1\) is retained for an ordinary output
scalar. At time seven the embedded slot supplies its color and the other
nine large slots supply their prices. For example, at a time other than
zero or seven the blocking row \(2302\) has weight \(4m^7o^3/f^{72}\).
Its weight at time zero
is \(2m^7o^3/f^{72}\) with \(a=1\), and at time seven it is
\(4m^6o^3/f^{65}\) with \(b=1\).

For each of the one-row tables below, group eligible rows by the stated
key \(\kappa\). Store the integer
\begin{equation}\label{appcert:bucket}
 A_d(\kappa)=\mathcal R\!\left(S\sum_{r\text{ with key }\kappa} W_d(r)\right)
\end{equation}
and use zero for a missing key. Thus the sum inside \eqref{appcert:bucket}
is rational and all subsequent path multiplications and additions are
exact integer operations.

\paragraph{The two-row head and the base family.}
At \(d=0,1\) use the key \((i_r,j_r,x(r))\), retaining the blocking symbol.
For a pair of types put \(\widehat h(x,y)=h(x,y)\) when both types are
numerical and \(\widehat h(x,y)=\bot\) otherwise. The two-step head table is
\[
 H_a(i,j,h)=
 \sum_{\substack{p,x,y\\\widehat h(x,y)=h}}
 A_0(i,p,x) A_1(p,j,y).
\]
At \(d=2,\ldots,13\) use the key \((p,q,e)=(i_r,j_r,e(r))\).
For each suffix start \(j\), initialize \(V_2(j,0)=1\) and the other
entries to zero. At each time \(d\) start a fresh zero table and add
\[
 V_d(p,g) A_d(p,q,e)\quad\text{to}\quad V_{d+1}(q,g\oplus e),
 \qquad
 g\oplus e=\begin{cases}
 \bot &\text{if }g=\bot\text{ or }e=\bot,\\
 g+e &\text{otherwise}.
 \end{cases}
\]
Write \(V_{14}^{j,b}\) for the terminal table, with embedded color \(b\)
imposed at time seven. The base contribution \(L_{ab}\) is the sum of
\(H_a(i,j,h)V_{14}^{j,b}(i,g)\) over \(i,j,h,g\) for which \(h=\bot\),
\(g=\bot\), or both indices are numerical and \(\ell_h\le g\le u_h\). Thus a blocking row in either part accepts
the walk. Each update appends one row at its source port, so induction on
\(d\) counts every suffix with the indicated port and score. The terminal
port \(i\) closes the walk, and the fixed time origin counts each word once.

\paragraph{The two extensions.}
For each state \(\sigma\) use only non-blocking rows. Form the four score
sets of a one-step first companion by putting \(e_c\) in
\(B_\sigma(r)_{i_cj_c}\) for each \((c,d')\in C_\sigma(x(r))\).
In the suffix use the finer key \((p,q,e,B)=(i_r,j_r,e(r),B_\sigma(r))\)
in \eqref{appcert:bucket} and denote its integer table by \(A_d^\sigma\).
For a fixed suffix start \(j\), initialize \(V_2^\sigma(j,0,I)=1\),
where the diagonal entries of \(I\) are \(\{0\}\) and the off-diagonal
entries are empty. Initialize the other entries to zero. At each step,
in a fresh table, add
\[
 V_d^\sigma(p,g,M) A_d^\sigma(p,q,e,B)
 \quad\text{to}\quad V_{d+1}^\sigma(q,g+e,M*B),
 \qquad
 (M*B)_{\alpha\gamma}
     =\bigcup_{\nu=0}^1\{u+v:u\in M_{\alpha\nu},\ v\in B_{\nu\gamma}\}.
\]
Here \(g\) is the score of the candidate suffix, while
\(M_{\alpha\gamma}\) records attainable first-companion scores from
\(\alpha\) to \(\gamma\). Each set can be stored as a bit mask, with bit
\(w\) representing score \(w\).

Take \(g>u_h\) for \(\sigma=1\) and \(g<\ell_h\) for \(\sigma=2\),
always with \(0\le g\le24\). Put \(D_1=48,\ D_2=24\). For each
\((h',h'',\alpha,\gamma)\in\mathcal H_\sigma(h)\), test whether there is
a \(w\in M_{\alpha\gamma}\) with \(w'=D_\sigma-g-w\in[0,24]\) and
\[
  \begin{cases}
  w\le u_{h'},\ w'\le u_{h''} &(\sigma=1),\\
  w\ge\ell_{h'},\ w'\ge\ell_{h''} &(\sigma=2).
  \end{cases}
\]
If such a \(w\) exists for any head pair, discard the terminal state for
this \(h\). Every genuine compatible pair obeying the two bounds supplies
a witness to this test. For numerical \(h\), define \(E^\sigma_{h,ab}\)
as the sum of \(H_a(i,j,h)V_{14}^{\sigma,j,b}(i,g,M)\) over
\(i,j,g,M\) in the indicated outside range that survive the test.
These sums give the weight of each candidate extension before selecting
its heads.

\paragraph{Selecting the heads and forming an entry.}
Choose \(I_1,I_2\) using the main pass, separately for each state.
Give head \(h\) weight \(w_h=\sum_{a,b} E^\sigma_{h,ab}\). Join the heads
at each pair of distinct positions of every compatible head triple,
giving a loop to a head that occurs twice. Start with the set \(U\) of
loopless heads of positive weight. While \(U\) is nonempty, select the
head \(h\in U\) that maximizes, in lexicographic order,
\[
 \left(\frac{w_h}{1+\sum_{j\in U\cap N(h)}w_j},\,w_h,\,h\right),
\]
where \(N(h)\) is the set of neighbors of \(h\). Remove \(h\) and its
neighbors from \(U\). The selected heads form \(I_\sigma\). Each
compatible triple has at most one selected occurrence, as is also
checked directly on the 468 triples. Both tilted passes use the same
\(I_1,I_2\) chosen by the main pass in that round.

For each pass, using its own row weights and rounding map, form the entry
\begin{equation}\label{appcert:entry}
 P_{ab}=
 \mathcal R\!\left(\frac{f L_{ab}}{S^{13}}\right)
 +\sum_{\sigma=1}^2 \mathcal R\!\left(
   \frac{f}{S^{13}}\sum_{h\in I_\sigma}E^\sigma_{h,ab}\right).
\end{equation}
The base family and the upper and lower extensions are disjoint.
Each of their contributions is converted separately to the output scale.
A fourteen-row product carries scale \(S^{14}\) and denominator
\(f^{13\cdot72+65}=f^{1001}\). The factor \(f/S^{13}\) leaves scale
\(S\) and accounts for the deleted scalar, giving normalization by
\(f^{1000}\). For fixed \(I_1,I_2\), each family sum is a polynomial with
nonnegative integer coefficients in its one-row table entries.
Consequently the rounding maps in \eqref{appcert:bucket} and \eqref{appcert:entry}
give the following entrywise bounds.
\begin{equation}\label{appcert:entrywise}
 P_0\ \le\ \frac{S}{f^{1000}}K_0,\qquad
 P_1\ \ge\ \frac{S}{f^{1000}}K_-,\qquad
 P_2\ \ge\ \frac{S}{f^{1000}}K_+ .
\end{equation}

\subsection{Integer spectral comparison and the final round}\label{appcert:final}

For \(P_0=\left(\begin{smallmatrix}a&b\\c&d\end{smallmatrix}\right)\), set
\[
 J=\left\lfloor
       \frac{a+d+\lfloor\sqrt{(a-d)^2+4bc}\rfloor}{2}
    \right\rfloor=\lfloor\rho(P_0)\rfloor .
\]
For a tail matrix \(Q=\left(\begin{smallmatrix}u&v\\w&x\end{smallmatrix}\right)\),
the two integer tests \(J>\max(u,x)\) and \((J-u)(J-x)>vw\) put \(J\)
strictly above its larger eigenvalue. Applying these tests to both
tails, entrywise monotonicity of the Perron root and
\eqref{appcert:entrywise} give
\[
 \max\{\rho(K_-),\rho(K_+)\}
 < (f^{1000}/S)J \le \rho(K_0).
\]
Theorem~\ref{thm:spectral} then yields \(\Lambda\ge f(J/S)^{1/1000}\).
For an integer \(u\), the two comparisons
\begin{equation}\label{eq:cert-floor}
 S(2500u)^{1000}<J(5573G)^{1000}<S\bigl(2500(u+1)\bigr)^{1000}
\end{equation}
place this certified endpoint strictly between \(u/G\) and \((u+1)/G\).
In particular \(u/G<\Lambda\) can be used as the seed of another round.

In Round~9 both selected sets have 78 heads. With rows indexed by
output color and columns by embedded color, each in the order \((O,T)\),
the three integer matrices from \eqref{appcert:entry} are
\[
\begin{aligned}
 P_0&=\begin{pmatrix}
 3883434916318410976492736782 & 5941763549053755133445913374\\
 4016360405175069604838470943 & 6145925632713824607452932023
 \end{pmatrix},\\[4pt]
 P_1&=\begin{pmatrix}
 3675111538177422381906390254 & 5624569415153833293664040644\\
 3796727775758173753987808316 & 5811434909041351886978678470
 \end{pmatrix},\\[4pt]
 P_2&=\begin{pmatrix}
 3863371480567745730128681998 & 5910846959000521536454774218\\
 3996202103131144563146837979 & 6114852651846792501758843116
 \end{pmatrix}.
\end{aligned}
\]
Applying the integer square-root formula to the three displayed matrices
gives, in the same order,
\[
\begin{aligned}
 \lfloor\rho(P_0)\rfloor&=10029057595358243801843755621=J,\\
 \lfloor\rho(P_1)\rfloor&=\phantom{0}9486260500781793939302943942,\\
 \lfloor\rho(P_2)\rfloor&=\phantom{0}9977922637375534216990157351.
\end{aligned}
\]
Each tail radius is strictly smaller than its displayed floor plus one.
Dividing by \(S\) gives the convenient upper bounds \(0.94862606\) and
\(0.99779227\) for the two normalized tail radii, while the normalized
main radius is at least \(J/S=1.0029057595358243801843755621\).
Using this \(J\) in \eqref{eq:cert-floor} gives
\(\Lambda>2.229206468135675727981404\).

\subsection{Iteration Table and Verification Commands}

Table~\ref{tab:cert-chain} lists $J$ and $u/G$ for Rounds~$0$ through $9$. Round~$0$ exceeds $11/5=2.20$, so Round~$1$ uses $2.20$ for ease of independent replay; from Round~$2$ onward, each round takes the full $24$-decimal floor from the previous row.

\begin{table}[!htbp]
\centering
\caption{Per-round integer Perron root floor $J$ and certified capacity lower bound $u/G$}
\label{tab:cert-chain}
{\small
\begin{tabular}{@{}r@{\qquad}r@{\qquad}l@{}}
\toprule
Round & Root floor $J$ & Certified lower bound $u/G$ \\
\midrule
$0$ & $8619997798499300956166357653$ & $2.228868987790439509571012$ \\
$1$ & $9845319270266868168558992989$ & $2.229165249379640409502069$ \\
$2$ & $10028797434351935063542336616$ & $2.229206410307698641566500$ \\
$3$ & $10029057230362894630468021680$ & $2.229206468054546469298121$ \\
$4$ & $10029057594846176740447243000$ & $2.229206468135561908393216$ \\
$5$ & $10029057595357525401811520045$ & $2.229206468135675568299202$ \\
$6$ & $10029057595358242793970674297$ & $2.229206468135675727757379$ \\
$7$ & $10029057595358243800429767335$ & $2.229206468135675727981090$ \\
$8$ & $10029057595358243801841773747$ & $2.229206468135675727981404$ \\
$9$ & $10029057595358243801843755621$ & $2.229206468135675727981404$ \\
\bottomrule
\end{tabular}}
\end{table}

Round~8 already reaches the 24-decimal lower bound of
Theorem~\ref{thm:main}. Taking that value as the seed of Round~9
gives the matrices displayed above. In the electronic supplement,
\path{verify_appendix_a.py} implements the recurrences of this appendix
using the Python standard library. With the script and
\path{appendix_a_data.json} in the same directory, the commands
\begin{quote}\ttfamily
python -B verify\_appendix\_a.py --round 9\\
python -B verify\_appendix\_a.py
\end{quote}
check the final round at its stated seed and the full chain from \(z=2\),
respectively. The full replay checks the compatible triples, the exact
price tests, the two tail comparisons and both strict inequalities
in \eqref{eq:cert-floor} at every round.

\clearpage
\section{Seven-Dimensional Point Table and Pairwise Checks}\label{app:seven-alphabet}

Encode $v=(v_0,\ldots,v_6)\in\F_3^7$ as $[v]=\sum_{j=0}^6 v_j 3^j \in \{0,\ldots,2186\}$. Projecting the $8$-dimensional $512$-point cap set \cite{funsearch} along coordinate $0$ (using the coordinate order of \href{https://github.com/google-deepmind/funsearch/blob/main/cap_set/n8_size512.txt}{\texttt{n8\_size512.txt}}, archived at \path{inner/n8_size512.txt}) gives $160$ two-point fibers (forming $T_0$, where preimages have $0$-th coordinate $1,2$) and $192$ single-point fibers (coordinate $0$). In Table~\ref{tab:T0-hex}, bit $k$ (LSB $=0$) of the hex integer in row $j$, column $i$ is $1$ iff $243j+81i+k \in T_0$ ($0\le i\le 2,\ 0\le k<81$).

\begin{table}[!htbp]
\centering
\caption{Hexadecimal bitmask encoding of the $160$-point seed set $T_0$}
\label{tab:T0-hex}
{\small
\begin{tabular}{@{}rlll@{}}
\toprule
$j$ & $i=0$ & $i=1$ & $i=2$ \\
\midrule
0 & \texttt{6c36000} & \texttt{4e000009c000000001b0}  & \texttt{27000004e0000001b0} \\
1 & \texttt{1008180} & \texttt{1000000200000000040}   & \texttt{8000001000000040} \\
2 & \texttt{201030}  & \texttt{1000000200000000008}   & \texttt{8000001000000008} \\
3 & \texttt{100920}  & \texttt{1000000200000000004}   & \texttt{8000001000000004} \\
4 & \texttt{40244}   & \texttt{105800020b0000000001}  & \texttt{b0080016010000000001} \\
5 & \texttt{4020c}   & \texttt{80c80010190000000001}  & \texttt{190400032080000000001} \\
6 & \texttt{80490}   & \texttt{1000000200000000002}   & \texttt{8000001000000002} \\
7 & \texttt{40242}   & \texttt{20980004130000000001}  & \texttt{130100026020000000001} \\
8 & \texttt{4020a}   & \texttt{100d000201a0000000001} & \texttt{1a0800034100000000001} \\
\bottomrule
\end{tabular}}
\end{table}

The four auxiliary subsets in $M = T_0 \cup E \cup C \cup D_+ \cup D_-$ are:
\[
\begin{aligned}
  E &= \{9, 18, 81, 162\},\\
  D_+ &= \{4, 972, 1952\}, \qquad D_- = \{8, 976, 1944\},\\
  C &= \{247, 494, 733, 973, 975, 980, 1222, 1466, 1706, 1946, 1948, 1950\}.
\end{aligned}
\]
Coloring and verification proceed as follows:
\begin{enumerate}
  \item \textbf{Ordinary set $U$:} $U$ has $200$ points, $80$ proper lines, and $64$ isolated points. Among the remaining $136$ points, $8$ have degree $10$:
  \[
    1007,\ 1034,\ 1247,\ 1274,\ 1735,\ 1762,\ 1975,\ 2002.
  \]
  Assigning color $0$ to these $8$ points and connecting the two uncolored endpoints of each proper line produces a bipartite graph with two $64$-point parts colored $+1$ and $-1$ (choosing $+1$ for the smallest index in each component), so every line has all three colors.

  \item \textbf{Marked set $M$:} Assigning $0$ to $E\cup C$ and $+1,-1$ to $D_+,D_-$, and propagating the zero-sum line rule with $T_0\subset\{\pm 1\}$ across components yields $258$ three-colored lines and $8,35,35$ internal lines in the zero, positive, and negative fibers. Each $83$-point nonzero fiber has $24$ isolated points; propagating the same rule on the remaining $59$ points (with $4,1952,972$ colored $-1,0,+1$ on the positive side and $8,976,1944$ on the negative side) yields second-level fibers of sizes $29,1,29$. Of the $35$ lines on each side, $13$ are three-colored and the remaining $22$ split evenly ($11$ each) inside the two $29$-point fibers, each forming an $11$-spoke star plus $6$ isolated points.

  \item \textbf{Zero-fiber projection:} All $8$ lines of $E\cup C$ lie in $C$, with every point of $C$ on two lines. Checking $f(v)=v_1+v_6$ confirms $\{f(a),f(b),f(c)\}=\{0,1,2\}$ on every line and $4$ preimages per value in $C$.
\end{enumerate}
The script \path{verify_geometry_and_flags.py} verifies the projection and calls \path{inner/verify_seven_leaf_sections.py} and \path{outer/verify_six_physical_marks.py} over all triples $[a]<[b]<[-a-b]$.

\clearpage
\section{Binary flags and acceptance probabilities of grouped nodes}\label{app:flag-node}

The coloring in Lemma~\ref{lem:physical-node} has two numerical inputs:
the probability of accepting a word with one specified zero and the probability
of accepting a word with only nonzero coordinates. We construct the coloring,
prove its recovery property, and derive both probabilities for the same choice
of layer orientations. The finite inputs are the four pivot matchings in
\path{matching.tsv} and the reader choices in \path{reader_rows.tsv}, both in
the electronic supplement \path{appendix_c_data.zip}.

\subsection{From a finite pivot matching to all layers}\label{appflag:matching}
Write \([b]=\{1,\ldots,b\}\). A word with one zero at \(i\) and plus set
\(S\) corresponds to the edge \(S\to S\cup\{i\}\) of the Boolean lattice.
All other entries equal \(-1\). Draw one uniform \(V\in\mathrm{GL}(b,2)\).
For a set \(I\) of column indices define
\[
 P_V(I)=\{j\in[b]:\operatorname{rank}V_{\le j,I}>
                          \operatorname{rank}V_{\le j-1,I}\}.
\]
The rank of the empty prefix is zero. Fix an orientation for each layer
\(h=|S|\). In the forward orientation use the displayed word and matrix.
In the reverse orientation first interchange its plus and minus signs, then
use \(J(V^{-1})^\top\), where \(J\) reverses row order. Both matrices are uniform.

Write \(T\) for the set of plus coordinates after orientation, \(s=|T|\),
\(P=P(T)\), and \(Q=P(T\cup\{i\})\), using the matrix for that orientation.
Restricting a row dependence from \(T\cup\{i\}\) to \(T\) gives
\(P\subset Q\). Invertibility gives \(|P|=s\) and \(|Q|=s+1\).
Thus \(Q\setminus P\) consists of a single row.

\paragraph{The pivot matching.}
At \(s=0\) accept every one-zero word. These words use only \(0\) and
\(-1\) in the oriented alphabet, so they themselves form a cap.
At \(s=1\), write \(P=\{p\}\) and accept precisely when
\(Q=\{p,p+1\}\). This is a matching. A step from \(\{p\}\) through its
matched upper end to the other singleton strictly increases \(p\).
In the model of independent fair binary rows, the sole old pivot has
distribution \(\Pr(p=j)=2^{-j}\). The new column passes the \(j-1\)
earlier nonpivot rows and first leaves the old span at \(j+1\) with
conditional probability \(2^{-j}\). The acceptance probability is
therefore \(\sum_{j\ge1}4^{-j}=1/3\).

For \(s\ge2\), choose the window parameters from the following table.
The supplement \path{matching.tsv} gives the partial function
\(\tau_{m,w}\) used in each row.
\[
\begin{array}{c|rrrr}
 s                 & 2 & 3,4 & 5,6 & \ge7\\
 \hline
 m                 & 2 & 3 & 5 & 7\\
 w                 & 6 & 7 & 10 & 15\\
 |\operatorname{dom}(\tau_{m,w})| & 22 & 61 & 586 & 15\,828
\end{array}
\]
Freeze the first \(s-m\) pivots of \(P\), and let \(u\) be the last
frozen pivot, with \(u=0\) for an empty prefix. Form the set of relative
row indices
\[
 H=\{j\in[w]:u+j\in P\},\qquad
 M(H)=\sum_{j\in H}2^{j-1}.
\]
Each entry gives \(M(H)\) in hexadecimal and the target
\(t=\tau_{m,w}(H)\in[w]\setminus H\) in one-based row numbering.
Accept precisely when the entry exists, \(u+t\le b\), and
\(Q=P\cup\{u+t\}\). The parameter \(w=15\) counts consecutive row
positions following the prefix. Any of the seven remaining pivots beyond
those positions still belong to \(P\).
For example, at \(s=7\) and \(b\ge20\), the profile
\(P=\{1,2,16,17,18,19,20\}\) has \(u=0\), \(M(H)=3\), and
\(\tau_{7,15}(\{1,2\})=3\). This edge is selected exactly when the
added pivot is row \(3\).

\paragraph{Verification on one window.}
For each listed \(H\) put \(U(H)=H\cup\{\tau_{m,w}(H)\}\).
Form a directed graph on the listed sets of each cardinality by drawing
\(H\to K\) whenever \(K\ne H\) and \(K\) is a facet of \(U(H)\)
(a subset obtained by deleting one element).
The following finite check uses all entries of each of the four tables.
\begin{verbatim}
seen = empty set
for each listed pair (H,t):
    check H occurs once, H is contained in [w], and |H| <= m
    check t belongs to [w] and t does not belong to H
    U = H union {t}
    check U is absent from seen, then insert U into seen
    N[H] = all listed facets of U other than H
for each H and each K in N[H]:
    for each L in N[K]:
        check H is absent from N[L]
\end{verbatim}
The four graphs have respectively \(36,134,2052,89\,288\) directed
edges, all with distinct ends. The checks give distinct upper ends and
no directed triangle. An alternating six-cycle would give a directed
triangle by traversing each selected edge upwards and the following
incidence downwards.

\paragraph{Lifting the window check.}
The upper end \(P\cup\{u+t\}\) has the same first \(s-m\) entries as \(P\),
since \(u+t>u\). It therefore determines the prefix, the window, and
all pivots outside the window. Its intersection with the window
determines \(H\) by uniqueness of \(U(H)\). Thus the full upper end
determines its unique matched lower end.
Now traverse a matched edge upwards and descend to a different
\(s\)-element facet. Deleting an entry from the upper end can only
increase the first \(s-m\) order statistics, coordinate by coordinate.
On a closed walk each is constant. Hence the same window is used
throughout the walk. Every upward step adds a pivot inside that window,
so the set of pivots beyond \(u+w\) can only shrink along the walk.
It, too, is constant on a closed walk. With both the prefix and the
pivots beyond the window fixed, distinct lower ends give distinct
window sets \(H\). A full alternating six-cycle would now project to
a directed triangle in the checked graph. This proves the asserted
matching and six-cycle exclusion in every positive oriented rank.

\subsection{One-zero words of a fixed coordinate sum}\label{appflag:onezerocap}
Fix \(\pi(w)=\sum_i w_i\in\F_3\). If two words in a zero-sum triple
have their zero at the same coordinate, the third has its zero there too.
At all remaining coordinates the three signs must agree. In a proper
line of one-zero words the zeros therefore occupy three distinct
positions \(i,j,k\); outside these positions the words agree.
The number \(h\) of pluses satisfies \(\pi(w)=2h-(b-1)\pmod 3\).
The three counts differ by at most two, hence are equal and the three
words use the same orientation. After that orientation each of the
three positions contains a plus, a minus and a zero, and each word
contains one of each there. If \(A\) is the common plus set outside
these positions, the selected edges must have the cyclic form
\[
 A+j\ \to\ A+i+j,\qquad A+k\ \to\ A+j+k,\qquad A+i\ \to\ A+k+i ,
\]
where \(A+i=A\cup\{i\}\). The other three incidences between
these lower and upper sets complete a six-cycle.

The pivot map preserves these incidences. If the three image edges
are distinct, the matching makes their three lower and three upper
ends distinct, giving the forbidden alternating six-cycle.
If exactly two image edges coincide, the two different lower pivot
sets are subsets of two different upper pivot sets. Their union,
of size at least \(s+1\), would have to equal both upper sets.
It remains to consider three identical image edges with added pivot \(p\).
For \(t=p-1,p\), form the quotient \(E_t\) of the first \(t\)
row coordinates by the span of the truncated columns of \(A\).
At \(p-1\), each of the three columns \(i,j,k\) lies modulo this
span in the line generated by the next. Their images are either all
zero or the same nonzero vector, since the field is binary.
Row deletion induces a surjection \(E_p\to E_{p-1}\) whose kernel
has dimension at most one. The three columns have at most two
possible lifts of their common image, so two have equal images in
\(E_p\). The edge joining this pair has gained no rank modulo its
lower span by row \(p\), a contradiction. This proves that the
accepted one-zero words in the fixed sum class form a cap. The
all-nonzero words are contained in \(\{\pm1\}^b\), itself a cap.

\subsection{Mixed lines and the choice of readers}\label{appflag:tables}

A mixed line of fixed coordinate sum has two one-zero words and one
all-nonzero word. With just one one-zero word, the other two have
opposite signs at its zero and agree elsewhere, giving three different
coordinate sums. In the remaining case the two zeros occupy different
coordinates, say \(i,j\). At those coordinates the all-nonzero word
has the same sign: its two differences from either one-zero word sum
to zero in \(\F_3\), which forces the two signs to agree. Write \(X\)
for its positive set and \(k=|X|\). If \(i,j\in X\), the two edges are
\[
 (X\setminus\{i,j\})\longrightarrow X\setminus\{i\},\quad
 (X\setminus\{i,j\})\longrightarrow X\setminus\{j\}.
\]
They belong to layer \(k-2\) and share their lower end. If
\(i,j\notin X\), the two lower ends are \(X\cup\{i\}\) and
\(X\cup\{j\}\); both edges end at \(X\cup\{i,j\}\) and belong to
layer \(k+1\). Sign reversal interchanges the two descriptions.

For a column set \(I\), restriction to the pivot rows \(P(I)\) is an
isomorphism on its span. Denote by \(v_p(I)\) the unique vector whose
restriction is the unit vector at \(p\in P(I)\). On a selected edge
\(S\to U=S\cup\{i\}\), let \(d\) be its added pivot. A binary
\emph{upper reader} \(f(U)\) takes component \(r\) of \(v_d(U)\).
A binary \emph{lower reader} \(g(S)\) takes component \(d\) of
\(v_r(S)\). The prescribed \(r\), or the list of rows to be XORed,
depends only on the layer, on the end in question, and on its pivot
set. The matching determines \(d\) from either end in positive oriented
rank. Give the edge color \(+1\) or \(-1\) according as
\(f(U)+g(S)\) in \(\F_2\) equals \(0\) or \(1\). Set unassigned
reader functions to zero, including those at oriented rank zero.
Give an all-nonzero word color zero and retain it if the colors on
every selected mixed pair through it in its coordinate-sum class differ.

These rules already imply the recovery assertion of
Lemma~\ref{lem:physical-node}. A mixed line that survives has the
three colors \(0,+1,-1\). Three points of a single input type with
the same coordinate sum cannot form a proper line by the preceding
subsection. Hence a zero-sum triple of accepted inputs in that
class has zero-sum outputs, and equal outputs force equal inputs.
It remains to count the inputs retained by one choice of the readers.

\paragraph{Which readers can be chosen independently.}
On a common-lower pair the two lower reads cancel; on a common-upper
pair the upper reads cancel. For fixed layer \(h\) and fixed orientation,
one reader is thus used only at \(k=h+2\) and the other only at \(k=h-1\).
This permits a separate choice at each \(k\) once the orientations of the
two incident layers are specified. We give no all-nonzero credit when
\(\min(k,b-k)<9\). In the range \(9\le\min(k,b-k)<16\), use just
the pivot set \([k]\). Its reversed complementary profile is
\([b-k]\). A deletion endpoint of the contiguous profile determines
its single missing row, and an addition endpoint determines its single
extra row. Use these requested rows as readers. The row computations
in Section~\ref{appflag:occupancy} show that they distinguish every
mixed pair. In particular these choices impose no conditions on the
echelon coefficients. Every one-zero endpoint used in either
orientation in these priced layers has lower rank at least seven.

\paragraph{Translating the central profile.}
For \(k,b-k\ge16\), put \(u=k-16\). For the all-nonzero words
we price the profiles
\begin{equation}\label{appflag:padding}
 P=[u]\ \cup\ (u+T),\qquad T=\{t_1<\cdots<t_{16}\}\subseteq[32],
 \qquad \delta(T)=\sum_{i=1}^{16}(t_i-i)\le30 .
\end{equation}
Here \([0]\) is empty. Subtract \(u\) to express a row of the
central block in local coordinates. If \(W\) is the span of the columns indexed by \(X\), the columns
indexed by \(X^c\) in \(J(V^{-1})^\top\) span \(JW^\perp\), since
the columns of \((V^{-1})^\top\) are dual to those of \(V\). A basis for the orthogonal
complement has identity coordinates on the original nonpivot rows;
after reversal these become its pivot rows. Its leading contiguous
block has length \(b-k-16\) and its central profile is
\(T^*=\{33-j:j\in[32]\setminus T\}\). In particular
\(\delta(T^*)=\delta(T)\). When both orientations are reversed
we use this latter profile as the starting profile.

Here is an explicit list of all requests made by a central profile.
For a local lower profile \(A\), let \(\tau(A)\) be the new row
specified by the seven-pivot window, measured in the same local
coordinates as \(A\); omit a request whenever the matching is
undefined.
\begin{enumerate}
\item For a deletion, take each \(e\in T\) and set \(Q=T\setminus\{e\}\).
If some \(d\in Q\), \(d<e\), has \(\tau(Q\setminus\{d\})=d\),
record an upper-end request with key \(Q\) and ideal read row \(e\).
There is at most one such \(d\), by upper-end uniqueness.
\item For an addition, take each \(a\notin T\) with
\(1\le a\le\max T+15\). Set \(Q=T\cup\{a\}\).
If \(d=\tau(Q)>a\), record a lower-end request with key \(Q\)
and ideal read row \(a\).
\end{enumerate}
The list applies equally to \(T^*\). To see that padding loses no
actual pairs, first consider a deletion of \(d<e\). Even after the
two deletions there are at least fourteen pivots above the
contiguous prefix. The seven-pivot window starts after the prefix,
so its target \(d\), and hence \(e\), must be above the prefix.
For an addition, the prefix has no vacant rows. The anchor for
\(\tau(T\cup\{a\})\) is at most \(\max T\), so \(a<d\) forces
\(a\le\max T+15\). In both cases deleting the contiguous prefix
translates the anchor, the window and its target by the same amount.
The depths measured downwards in pivot order or upwards in nonpivot
order are also preserved: earlier pivots do not precede a deletion
in reverse pivot order, and the prefix has no free rows.
This proves the completeness of the two finite loops.

There are \(26\,103\) sets \(T\) in \eqref{appflag:padding}; they can
be enumerated by taking nondecreasing integers \(t_i-i\) in
\([0,16]\) with sum at most \(30\). The electronic table
\path{reader_rows.tsv} gives the final choices of reads for their
requests. For equal orientations the two legs use deletion and
addition in a single flag. For forward then reverse they both use
deletion, and for reverse then forward they both use addition. For
the latter two tables leg 0 uses \(T\) and leg 1 uses \(T^*\).
The entries are indexed by the leg and by the hexadecimal integer
\(\sum_{j\in Q}2^{j-1}\); their row numbers are local and start at
one. Add back the prefix length in the flag of that leg to evaluate
the reader on the actual subspace. Numbers of keys on legs 0 and 1
respectively are \(33\,761,35\,457\) for the equal-orientation table,
\(33\,761,33\,761\) for two deletions and \(35\,457,35\,457\) for
two additions. For instance, key \texttt{17baff} of leg 1 in
the equal-orientation table reads rows \(12\) and \(19\) by XOR.
All other keys read one row. Applying the two loops above checks
every key directly.

\subsection{Reader mismatches and pairs of actual columns}\label{appflag:occupancy}

Fix an all-nonzero word with positive-coordinate set $X$, where $|X|=k$.
Write $W$ for the span of the corresponding columns of $V$, and $P=P_V(X)$ for its pivot rows.
Restriction to $P$ identifies $W$ with $\F_2^k$. For $p\in P$, the canonical vector
$v_p(X)$ restricts to $e_p$ on $P$. Its free entries are
$c_{pf}=v_p(X)_f$ for $f\notin P$ and $f>p$. Conditional on $P$, these bits are independent
and fair. Each assignment specifies one subspace with pivot set $P$, and all such
subspaces have the same number of ordered bases and completions to an invertible matrix.

\paragraph{The two readers.}
An upper-end request occurs when two one-zero words have the same lower end.
Deleting either of two columns of $X$ gives an upper end with pivot set $Q=P\setminus\{e\}$.
Deleting both gives $Q\setminus\{d\}$, where $d<e$ and the matching sends
$Q\setminus\{d\}$ to $Q$. The two upper-end canonical vectors indexed by $d$ sum to
$v_e(X)$. Indeed, their sum vanishes on $Q$ and is nonzero:
the two upper spaces meet in the common lower space, whose pivots omit $d$.
On $P$ the sum must therefore be the unit vector at $e$. A reader taking the bit in row $r$
therefore distinguishes the endpoints precisely when $v_e(X)_r=1$. The requested row $e$
always works. For another prescribed row this entry is either zero or a free bit $c_{er}$.

For a lower-end request the two words have the same upper end. Adding either of two
columns outside $X$ produces the same lower pivot set $Q=P\cup\{a\}$. Adding both
introduces one further pivot $d>a$, determined by the matching at $Q$. Let $w_i,w_j$
be the two added columns after their components on $P$ have been eliminated. They both
have first one at $a$, and $w_i+w_j$ has first one at $d$. In $W+\langle w_i\rangle$,
the canonical vector indexed by $r\in P$ equals $v_r(X)+v_r(X)_a w_i$,
and the vector indexed by $a$ equals $w_i$. Reading component $d$ in the two lower ends
thus gives difference $v_r(X)_a$ when the reader selects index $r\in P$, and difference
$1$ for $r=a$. If two reads are prescribed, their differences add in $\F_2$.

\paragraph{The affine equations.}
For each orientation mode, the reader table is indexed by the leg and the endpoint pivot
set $Q$. For each potential request impose that the difference of the two reads
equal one. In mixed orientations, the second leg uses the reversed orthogonal complement.
A free entry $(p,f)$ in its local $32$-row window is the entry $(33-f,33-p)$ in the original
window: the orthogonal-complement basis has unit coordinates on the original nonpivots,
and its values on the original pivots are given by the transposed free entries.
Reversing the rows gives the displayed indices. Repeated occurrences of a free entry
therefore use a single variable, also when the two requests come from different legs.

Write the resulting equations as $Ac=b_0$ over $\F_2$ and let $\sigma$ mark exactly the
failed equations. Then
\begin{equation}\label{appflag:syndrome-mass}
 \Pr(Ac+b_0=\sigma\mid P)=
 \begin{cases}
  2^{-\operatorname{rank} A},& b_0+\sigma\in\operatorname{im}A,\\
  0,&\text{otherwise}.
 \end{cases}
\end{equation}
Identically satisfied tests can be omitted. The certificate counts the zero syndrome and
the syndromes with one or two failed keys. These events are disjoint. For each nonzero
syndrome we require, in addition, that every failed key have no pair of actual columns
that could use it.

\paragraph{The forbidden pair for a key.}
For an upper-end key, let $\lambda_i$ be the linear form on $W$ taking value one on the
$i$th positive column and zero on the other positive columns. These forms are the dual
of the actual ordered basis of $W$. In pivot coordinates, deleting column $i$ restricts
$W$ to $\ker\lambda_i$. The removed pivot is the last nonzero entry of $\lambda_i$:
on that row the coordinate is determined by the preceding coordinates. For two deletions
with the same removed pivot $e$, eliminating row $e$ in the two equations shows that the
second removed pivot $d$ is the last nonzero entry of the sum of the two forms.
In \emph{decreasing} order on $P$, the two forms therefore have their first one at the
same position $j$ (that of $e$), and their first difference at position $\ell>j$
(that of $d$). For a lower-end key the columns outside $X$, reduced on $P$, form an
ordered basis of $\F_2^b/W$. In \emph{increasing} order on the nonpivots, the first
one of both reduced columns is at position $j$ (that of $a$), and their first difference
is at position $\ell$ (that of $d$). For either end it suffices to exclude a pair of this
type in the relevant basis. Within a leg distinct failed keys have distinct $j$:
for the fixed profile the removed or added row determines $Q$.

\paragraph{One or two failures on a single leg.}
For a basis of dimension $n$, let $\mathcal S$ be the set of depth pairs $(j,\ell)$
on this leg. Discard any with $\ell>n$, since they cannot occur. For $\mathcal S=\varnothing$
put $L_n(\mathcal S)=1$. Otherwise take $H=\max_{(j,\ell)\in\mathcal S}\ell$.
First sample $n$ independent uniform $H$-bit strings with replacement. Strings whose
first one is at $j$ fall into $m_{j\ell}=2^{\ell-j-1}$ bins according to the intervening
bits before $\ell$. Each bin contains two classes indexed by bit $\ell$, each with
$s_\ell=2^{H-\ell}$ patterns. A bin may use either class, but using both creates the
forbidden pair. The strings for different first-one positions are disjoint. Thus the
number of safe ordered $n$-tuples is
\begin{equation}\label{appflag:prefix-count}
 G_n(\mathcal S)=n![z^n]\exp(D_{\rm out}z)
       \prod_{(j,\ell)\in\mathcal S}(2\exp(s_\ell z)-1)^{m_{j\ell}},
 \qquad
 D_{\rm out}=2^H-\sum_{(j,\ell)\in\mathcal S}2^{H-j}.
\end{equation}
The factor for one bin counts assignments of labelled columns to either class and counts the empty
assignment once. For integer evaluation, expand the powers by the binomial theorem. Each choice
of integers $0\le v_{j\ell}\le m_{j\ell}$ contributes
\[
 \left(\prod_{(j,\ell)\in\mathcal S}
       (-1)^{m_{j\ell}-v_{j\ell}}{m_{j\ell}\choose v_{j\ell}}2^{v_{j\ell}}\right)
 \left(D_{\rm out}+\sum_{(j,\ell)\in\mathcal S}v_{j\ell}s_\ell\right)^n .
\]

The first $H$ rows of a uniform ordered basis have the law of this $H\times n$ matrix
conditional on full row rank: every rank-$H$ matrix has the same number of completions
by $n-H$ rows. There are $R_{H,n}=\prod_{i=0}^{H-1}(2^n-2^i)$ such matrices.
At most $2^{Hn}-R_{H,n}$ safe matrices can have deficient rank. Consequently
\begin{equation}\label{appflag:conditioned-prefix}
 L_n(\mathcal S)=
 \max\left\{0,\frac{G_n(\mathcal S)-2^{Hn}+R_{H,n}}{R_{H,n}}\right\}
 \ \le\ \Pr(\text{no requested pair in a uniform ordered basis}).
\end{equation}
For two keys $u,v$ on this leg use the larger of $L_n(\{u,v\})$ and
$L_n(\{u\})+L_n(\{v\})-1$.

For example, take $n=H=3$ and $u=(1,3)$. Formula~\eqref{appflag:prefix-count}
gives $4\cdot6^3-4\cdot5^3+4^3=428$ safe tuples out of $512$.
Adding $v=(2,3)$ gives $8\cdot5^3-12\cdot4^3+6\cdot3^3-2^3=386$.
There are $168$ full-rank matrices. The two lower bounds from
\eqref{appflag:conditioned-prefix} are $1/2$ and $1/4$. Enumerating the $168$
matrices gives $120$ and $96$ safe cases, respectively.

\paragraph{Failures on different legs.}
Conditional on the entire subspace $W$, the positive columns form a uniform ordered basis
of $W$. The images of the complementary columns form an independent uniform ordered basis
of $\F_2^b/W$: for each pair of bases there are $2^{k(b-k)}$ choices of lifts of the
latter basis. Conditional on $W$, the probabilities in pivot coordinates depend only on
the two dimensions and the specified depths. They can therefore be used after imposing
the syndrome in \eqref{appflag:syndrome-mass}, which depends only on $W$ for the fixed profile.
In $J(V^{-1})^\top$ the columns labelled by the original negatives give the reversed
dual basis to the quotient basis. The remaining dual columns, restricted to $W$ before
row reversal, give the basis dual to the original positive columns. These changes preserve
the two uniform laws and their conditional independence.

Writing $0,1$ for the forward and reverse orientations of the layers $k-2,k+1$, the
leg order of the certificate is as follows. Deletion and addition refer to columns in
the indicated flag.
\[
\begin{array}{c|ll}
 (e_{k-2},e_{k+1})&\text{first leg (dimension)}&\text{second leg (dimension)}\\ \hline
 (0,0)&\text{forward deletion }(k)&\text{forward addition }(b-k)\\
 (1,1)&\text{reverse deletion }(b-k)&\text{reverse addition }(k)\\
 (0,1)&\text{forward deletion }(k)&\text{reverse deletion }(b-k)\\
 (1,0)&\text{forward addition }(b-k)&\text{reverse addition }(k)
\end{array}
\]
For $(1,1)$ the profile and both reader equations are taken in the reversed flag, using the original negative-coordinate set of size $b-k$ as its positive-coordinate set. In the mixed cases only the second leg is folded back into the original window by the index map above.
A syndrome with one failed key on each leg thus receives the product of the two
single-key lower bounds in their indicated dimensions.

For a given profile let $p_P(\sigma)$ denote the resulting lower bound on the
conditional probability of avoiding all failed keys, with $p_P(0)=1$. The all-nonzero
contribution from this profile is bounded below by
\[
 \Pr(P)\!
 \sum_{\substack{|\operatorname{supp}\sigma|\le 2\\
                  b_0+\sigma\in\operatorname{im}A}}
           2^{-\operatorname{rank}A}\,p_P(\sigma).
\]
The summands refer to a fixed orientation pair. Distinct profiles are disjoint, so
their bounds add. Section~\ref{appflag:prices} gives the integer scaling used in this sum.

\subsection{Finite computation of both acceptance bounds}
\label{appflag:prices}

Fix a block length \(b\) and write \(D=2^{176}\). We first compute the acceptance bound for one specified zero coordinate when the chosen orientation has \(s\) positive columns. We then select all layer orientations jointly, using the same choices to price the all-nonzero words.

\paragraph{Pivot probabilities.}
For a pivot set \(P=\{p_1<\cdots<p_a\}\subseteq[a+e]\), put
\(\delta(P)=\sum_{i=1}^{a}(p_i-i)\).
The reduced column basis with these pivot rows has \(ae-\delta(P)\) free binary entries. All subspaces of dimension \(a\) are equally likely. Thus
\[
 \Pr(P)=\frac{2^{-\delta(P)}}{G(a,e)},\qquad
 G(a,e)=\sum_{\substack{P\subseteq[a+e]\\|P|=a}}2^{-\delta(P)}
       =\prod_{i=1}^{a}\frac{1-2^{-(e+i)}}{1-2^{-i}} .
\]
The displacements form a partition in an \(a\)-by-\(e\) rectangle. Separating partitions according to whether their first column is full gives
\(G(a,e)=G(a-1,e)+2^{-a}G(a,e-1)\).
The boundary values are \(1\), and \(G(a,e)=G(e,a)\). Let \(L(a,e),U(a,e)\) be integer lower and upper bounds for \(DG(a,e)\). Both boundary values equal \(D\). Swap the arguments if \(a>e\), then apply the recurrence to the two bounds, rounding the term divided by \(2^a\) down for \(L\) and up for \(U\). Conditional on the span of the old columns, reduction of the extra column on its pivot rows gives a uniform nonzero vector on the \(b-s\) free rows. The probability that its first nonzero entry is the \(v\)-th free row equals \(2^{-v}/(1-2^{-(b-s)})\).

\paragraph{One zero and at least seven positive columns.}
Read an entry of the seven-pivot matching as \((M,t)\), with \(M\subseteq[15]\), \(i=|M|\le7\), and the selected new pivot \(t\in[15]\setminus M\), using one-based row numbers. Write \(v=t-|\{p\in M:p<t\}|\). For a suffix of length \(m\), let \(W(m)\) be its integer weight. For \(7\le m<15\), sum \(D2^{-\delta(M)-v}\) over entries with \(i=7\) and \(M\cup\{t\}\subseteq[m]\). For \(m\ge15\), form
\[
 A_j=\sum_{\substack{(M,t)\\|M|=7-j}}
          D\,2^{-\delta(M)-v-j(15-|M|)},\qquad
 W(m)=\sum_{j=0}^{\min(7,m-15)}
      \left\lfloor\frac{A_j L(j,m-15-j)}{D}\right\rfloor .
\]
All summands defining \(A_j\) are integers. The last \(j\) pivots lie beyond row \(15\), where their displacements contribute the rectangular partition sum with width \(m-15-j\).

Put \(c=s-7\) and \(d=b-s\). When \(c>0\), write the last of the first \(c\) pivots as \(c+j\). The displacement of this last pivot contributes \(j\), the seven suffix pivots contribute an additional \(7j\), and the \(j\) free rows in the prefix must vanish in the reduced extra column. This gives the factor \(2^{-9j}\). The other \(c-1\) prefix pivots contribute \(G(c-1,j)\). Define the integer numerator
\begin{equation}
 T=\sum_{j=0}^{\min(d,21)}
 \left\lfloor 2^{-9j}
   \left\lfloor\frac{L(c-1,j)W(7+d-j)}{D}\right\rfloor\right\rfloor,
 \qquad
 a_b(s)=\frac{T}{U(s,d)}\frac{2^d}{2^d-1}.
 \label{appflag:onezero}
\end{equation}
If \(c=0\), take just \(j=0\) and replace \(L(c-1,0)\) by \(D\). The finite sum discards nonnegative terms. More explicitly, all rectangular sums are at most
\(\prod_{i\ge1}(1-2^{-i})^{-1}<4\), and the unrounded suffix weight divided by \(D\) is at most \(G(7,m-7)\). Hence terms with \(j\ge22\) contribute less than \(32\,2^{-198}/(1-2^{-9})\) to the unrounded conditional probability. The denominator above is rounded upwards.

\paragraph{Short ranks.}
For a window with \(r\) unresolved pivot columns and length \(n\), let \(\omega_r(M)\) be the probability of seeing the mask \(M\) in the first \(n\) independent fair binary rows. One obtains it by starting at rank deficit \(r\), multiplying by \(2^{-\ell}\) for a miss at deficit \(\ell\), and by \(1-2^{-\ell}\) for a pivot, after which the deficit decreases by one. The matching coefficient is
\[
 e_{r,n}=(1-2^{-(r+1)})
       \sum_{(M,t)}\omega_r(M)\,2^{-t+|\{p\in M:p<t\}|}.
\]
For \(s\ge r\), the probability in the infinite independent-row model for this window on the last \(r\) pivots is \(e_{r,n}/(1-2^{-(s+1)})\). Indeed, at deficit \(\ell\) each failed pivot row before the next pivot also requires one fair bit of the extra column to vanish. Averaging over the geometric waiting time gives \((1-2^{-\ell})/(1-2^{-(\ell+1)})\); the product over \(\ell=r+1,\ldots,s\) telescopes.

The finite matching tables give
\begin{align*}
 (e_{2,6},e_{3,7},e_{5,10})&=
 \left(\frac{10199}{32768},\frac{10510455}{33554432},
             \frac{2822461610100615}{9007199254740992}\right),\\
 e_* = e_{7,15}&=
 \frac{3254911033690990776061278018232125}
      {10384593717069655257060992658440192}.
\end{align*}
For \(2\le s<7\), the largest of the three applicable short-window rates (those with \(r\le s\)) is at least \(e_*\). The rank-zero and rank-one rates from Section~\ref{appflag:matching} also exceed \(e_*\). For the \(s+1\) specified columns the event of failure to have full rank in the first \(b\) rows has probability at most \(\varepsilon=(2^{s+1}-1)2^{-b}\), by the union bound over nonzero kernel vectors. Conditional on full rank these columns have exactly their law under uniform \(V\in\mathrm{GL}(b,2)\); the selected pivots have then been resolved within the first \(b\) rows. For \(s<7\) the computation therefore sets \(a_b(s)=(e_*-\varepsilon)/(1-\varepsilon)\).

\paragraph{Integer weights for the all-nonzero layers.}
Set \(M_0=2^{80}\). For a central profile \(T\) and a consistent
syndrome with at most two failed keys, Section~\ref{appflag:occupancy}
gives the unnormalized mass \(2^{-\delta(T)-\operatorname{rank}A}\).
Assign \(2^{80-\delta(T)-\operatorname{rank}A}\) to its bin \(F\), the
sorted list of triples (leg, first-one depth, first-difference depth).
Each exponent used is nonnegative. Write \(H_F\) for the sum in that
bin. In the order equal orientations, two deletions, two additions,
the numbers of nonempty bins are \(671,611,645\). As an intermediate
check, the respective values of \(H_{\varnothing}\) are
\[
\begin{array}{lr}
\text{equal orientations}&1370133050418933090222080\\
\text{two deletions}&1359607802008646048546816\\
\text{two additions}&1385606543351400544337920 .
\end{array}
\]

For a nonempty set of failures on one leg of dimension \(n\), put
\(Q_n(\mathcal S)=\lfloor D L_n(\mathcal S)\rfloor\) using
\eqref{appflag:conditioned-prefix}, and put \(Q_n(\varnothing)=D\).
The integer multiplier \(q_F\) is \(D\) for no failures and the
singleton \(Q_n\) for one failure. For two failures on one leg take the
larger of the joint \(Q_n\) and the sum of the singleton values minus
\(D\). For one failure on each leg, take the product of its two singleton
values, divide by \(D\), and round down. Use the dimensions indicated
in Section~\ref{appflag:occupancy}. For the orientation pair \(e,e'\)
the integer lower weight per word with \(k\) pluses is
\[
 C_k(e,e')=
 \begin{cases}
  0, &\min(k,b-k)<9,\\
  \lfloor D^2/U(k,b-k)\rfloor, &9\le\min(k,b-k)<16,\\
  \left\lfloor D\sum_F H_Fq_F/(M_0 U(k,b-k))\right\rfloor,
       &\min(k,b-k)\ge16.
 \end{cases}
\]
The second case prices the contiguous profile; the last prices the
central profiles in \eqref{appflag:padding}. Dividing each integer
weight by \(D\) gives a lower bound on its conditional probability.

\paragraph{Selecting the orientations.}
Give a one-zero vertex of height \(h\) and orientation \(e\) the integer weight
\[
 N_h(e)=\binom{b-1}{h}
          \left\lfloor D\,a_b\bigl((1-e)h+e(b-1-h)\bigr)\right\rfloor .
\]
An edge from \(h\) to \(h+3\) has weight \(\binom b{h+2}C_{h+2}(e,e')\). The boundary all-nonzero counts \(k=0,1,b-1,b\) have no such edge and receive zero credit. All one-zero heights \(0,\ldots,b-1\), including the ends, receive their vertex weights.

For each residue class modulo \(3\), start at its first height \(h\) with \(V_h(e)=(0,N_h(e))\). At successive heights take the lexicographic maximum
\[
 V_{h+3}(e')=\max_{e\in\{0,1\}}
       \left\{V_h(e)+\left(\binom b{h+2}C_{h+2}(e,e'),N_{h+3}(e')\right)\right\}.
\]
The two components count all-nonzero and one-zero mass, in that order. Maximize over the final orientation in each residue class, then add the three resulting pairs to obtain \((v,w)\); each maximizing path fixes a single actual choice of orientations on its residue class. Summing over independent fair signs gives
\[
        B_b=\frac{v}{2^b D},\qquad \eta_b=\frac{w}{2^{b-1}D}.
\]
Column-permutation invariance gives the same \(\eta_b\) for every specified zero coordinate. Flooring \(10^{12}\eta_b,10^{12}B_b\) for \(b=144,182,183,252,641\) gives \(313436531305\) in the first column each time and, in the same order,
\[
 328092514121,\quad328035264355,\quad328032524173,\quad327984670844,\quad327904358794
\]
in the second. The spectral calculation uses the two rational bounds before multiplication by \(10^{12}\) and flooring. The accompanying \path{verify_appendix_c.py} recomputes these ten
integers from the two rule tables by integer and rational arithmetic. It also compares all
five pairs of exact fractions with \path{node_rates.tsv}, the values used
by the spectral certificate.

\end{document}